\documentclass{article}
\usepackage{graphicx}
\usepackage{amsmath}
\usepackage{amsfonts}
\usepackage{amssymb}
\newtheorem{theorem}{Theorem}[subsection]

\newtheorem{claim}[theorem]{Note}

\newtheorem{conjecture}[theorem]{Conjecture}
\newtheorem{corollary}[theorem]{Corollary}

\newtheorem{definition}[theorem]{Definition}
\newtheorem{example}[theorem]{Example}

\newtheorem{lemma}[theorem]{Lemma}

\newtheorem{proposition}[theorem]{Proposition}

\newenvironment{proof}[1][Proof]{\textbf{#1.} }{\ \rule{0.5em}{0.5em}}

\begin{document}

\title{Kostka-Foulkes polynomials cyclage graphs and charge statistic for the root
system $C_{n}$}
\author{C\'{e}dric Lecouvey\\lecouvey@math.unicaen.fr}
\date{}
\maketitle
\begin{abstract}
We establish a Morris type recurrence formula for the root system $C_{n}%
$.\ Next we introduce cyclage graphs for \ the corresponding
Kashiwara-Nakashima's tableaux and use them to define a charge statistic.
Finally we conjecture that this charge may be used to compute the
Kostka-Foulkes polynomials \ for type $C_{n}.$
\end{abstract}

\section{Introduction}

There exists $q$-analogues of the multiplicities of weights in the irreducible
representations of the classical Lie algebras.\ They are obtained by
substituting the ordinary Kostant's partition function $\mathcal{P}$ by the
$q$-analogue $\mathcal{P}_{q}$ defined from the equality%
\[
\prod_{\alpha\text{ positive root}}\dfrac{1}{(1-qx^{\alpha})}=\sum_{\beta
}\mathcal{P}_{q}(\beta)x^{\beta}.
\]
Given $\lambda$ and $\mu$ two partitions%
\[
K_{\lambda,\mu}(q)=\sum_{\sigma\in W}(-1)^{l(\sigma)}\mathcal{P}_{q}%
(\sigma(\lambda+\rho)-(\mu+\rho)).
\]
As shown by Lusztig \cite{Lut} $K_{\lambda,\mu}(q)$ is a polynomial in $q$
with non negative integer coefficients. These polynomials naturally appear in
the classical theory of Hall-Littlewood polynomials.\ They coincide with the
Kostka-Foulkes polynomials that is, with the entries of the transition matrix
between the basis of Hall-Littelwood polynomials and the basis of Schur
functions \cite{mac}. Many interpretations of the Kostka-Foulkes polynomials
exist. For example, they appear in the filtrations of weight spaces by the
kernels of powers of a regular nilpotent element, and degree in harmonic
polynomials \cite{bry}, \cite{cal}, \cite{JLZ}. We recover them in the
expansion of the Hall Littlewood polynomials in terms of the affine Hecke
algebra (see \cite{NR}).\noindent

For type $A_{n-1}$ the positivity of the Kostka-Foulkes Polynomials can also
be proved by a purely combinatorial method. Recall that for any partitions
$\lambda$ and $\mu$ with $n$ parts the number of semi-standard tableaux of
shape $\lambda$ and weight $\mu$ is equal to the multiplicity of the weight
$\mu$ in the finite dimensional irreducible module of $U_{q}(sl_{n})$ with
highest weight $\lambda.$ In \cite{LSc1} Lascoux and Sch\"{u}tzenberger have
introduced a beautiful statistic $\mathrm{ch}_{A}$ on standard tableaux called
the charge and, by using Morris recurrence formula, have proved the equality%
\begin{equation}
K_{\lambda,\mu}(q)=\sum_{T\in ST(\mu)_{\lambda}}q^{\mathrm{ch}_{A}(T)}
\label{th_ls}%
\end{equation}
where $ST(\mu)_{\lambda}$ is the set of semi-standard tableaux of shape
$\lambda$ and weight $\mu.\;$Set $\mathcal{A}_{n}=\{1<\cdot\cdot\cdot
<n\}.\;$The charge may be defined by endowing $ST(\mu)$ the set of
semi-standard tableaux of weight $\mu$ with a structure of graph defined from
Lascoux-Sch\"{u}tzenberger's plactic monoid.\ Recall that the plactic monoid
is the quotient set of $\mathcal{A}_{n}^{\ast}$ the free monoid on
$\mathcal{A}_{n}$ by the Knuth relations%
\[
abx\equiv\left\{
\begin{tabular}
[c]{l}%
$bax$ if $a<x\leq b$\\
$axb$ if $x\leq a<b$%
\end{tabular}
\right.  .
\]
For any tableau $T$ we denote by $\mathrm{w}(T)$ the column reading of $T$
that is, the word obtained by reading the columns of $T$ from right to left
and from top to bottom. The cyclage graph structure on $ST(\mu)$ can be
defined as follows.\ We draw an arrow $T\rightarrow T^{\prime}$ between the
two tableaux $T$ and $T^{\prime}$ of $ST(\mu),$ if and only if there exists
$u$ in $\mathcal{A}_{n}^{\ast}$ and $x\neq1\in\mathcal{A}_{n}$ such that
$\mathrm{w}(T)\equiv xu$ and $\mathrm{w}(T^{\prime})\equiv ux$.\ Then we say
that $T^{\prime}$ is a cocyclage of $T.$ The essential tool to define this
graph structure is the insertion algorithm for the semi-standard tableaux. The
cyclage graph $ST(\mu)$ contains a unique row tableau $L_{\mu}$ which can not
be obtained as the cocyclage of another tableau of $ST(\mu)$. Let $T_{\mu}$ be
the unique semi-standard tableau of shape $\mu$ belonging to $ST(\mu)$. Then
there is no cocyclage of $T_{\mu}.$ For any $T\in ST(\mu)$ all the paths
joining $L_{\mu}$ to $T$ have the same length.\ This length is called the
cocharge of $T$ and denoted $\mathrm{coch}_{A}(T).$ Similarly, all the paths
joining $T$ to $T_{\mu}$ have the same length called the charge of $T.$ The
maximal value of $\mathrm{ch}_{A}$ is $\left\|  \mu\right\|  =\mathrm{ch}%
_{A}(L_{\lambda})=\sum_{i}(i-1)\mu_{i}.\;$Moreover the charge and the cocharge
satisfy the equality $\mathrm{ch}_{A}(T)=\left\|  \mu\right\|  -\mathrm{coch}%
_{A}(T)$ for any $T\in ST(\mu).$

\noindent The initial cocyclage of the tableau $T$ of reading $\mathrm{w}%
(T)=xu$ with $x\neq1$ is obtained by inserting $x$ in the sub-tableau of $T$
of reading $u.$ Every tableau $T\in ST(\mu)$ can be related to $T_{\mu}$ by a
sequence of initial cocyclages.\ So it is enough to consider initial
cocyclages to define $\mathrm{ch}_{A}\mathrm{..}$

\noindent The charge of $T$ can also be defined directly from $\mathrm{w}(T)$
when $\mu$ is a dominant weight. Moreover it can be characterized in terms of
the geometry of the crystal graph $B(\lambda)$ associated to $V(\lambda)$
\cite{LLT}.

In this article we restrict ourselves to the root system $C_{n}.\;$Our aim is
two folds.\ First we establish Morris type recurrence formula for type $C_{n}$
and use it to give explicit formulae for Kostka-Foulkes polynomials when
$\lambda$ is a row partition or a column partition of height $2.\;$Next we
introduce a cyclage graph structure and a notion of charge for type $C_{n}.$
For any dominant weight $\lambda,$ let $V(\lambda)$ be the finite dimensional
irreducible $U_{q}(sp_{2n})$-module with highest weight $\lambda$. In
\cite{KN} Kashiwara and Nakashima have given a combinatorial description of
$B(\lambda)$ the crystal graph of $V(\lambda)$ in terms of symplectic tableaux
analogous to the semi-standard tableaux for type $C_{n}$. From the plactic
monoid and the insertion algorithm described in \cite{Lec} it is natural to
try to obtain cyclage graphs for symplectic tableaux. Nevertheless the
situation is more complex than for type $A_{n-1}.\;$First we have to restrict
the possible cocyclage operations to the initial cocyclage to avoid loops in
our cyclage graphs. Moreover if we use the complete insertion algorithm for
type $C_{n},$ the number of boxes of the cocyclage of a tableau $T$ may be
strictly less than that of $T$ due to the contraction relation in the plactic
monoid. The cyclage graphs obtained by this mean seem to be not relevant to
define a charge related to the Kostka-Foulkes polynomials. To overcome this
problem we will execute the insertion algorithm without this contraction
relation and consider that the symplectic tableaux are filled by letters of
the totally ordered alphabet
\[
\mathcal{C}_{n}=\{\overline{n}<\cdot\cdot\cdot<\overline{1}<1<\cdot\cdot
\cdot<n\}
\]
which can be naturally embedded in the infinite alphabet%
\[
\mathcal{C}_{\infty}=\{\cdot\cdot\cdot<\overline{n}<\cdot\cdot\cdot
<\overline{1}<1<\cdot\cdot\cdot<n<\cdot\cdot\cdot\}.
\]
Our convention for the alphabet $\mathcal{C}_{n}$ is not identical to that of
\cite{KN} to dispose of a natural infinite extension of $\mathcal{C}_{n}%
$.\ Denote by $\mathbf{ST}(n)$ the set of symplectic tableaux defined on
$\mathcal{C}_{n}.\;$If $T\in\mathbf{ST}(n)$ the initial cocyclage (without
contraction) of $T$ does not belong to $\mathbf{ST}(n)$ in general but belongs
to $\mathbf{ST}(n+1).$ So it is natural to consider the cyclage graph
structure of $\mathbf{ST=}\underset{n\geq1}{\cup}\mathbf{ST}(n).$

\noindent Let $\mu$ be a dominant weight for the root system $C_{n}$.\ With
our convention $\mu$ may be identified with the partition $(\mu_{\overline{n}%
},...,\mu_{\overline{1}})$.\ We are going to endow $\mathbf{ST}(\mu)$ the
subset of $\mathbf{ST}$ containing the symplectic tableaux $T$ such that for
any $m\geq1,$ the number of letters $\overline{k}$ in $T$ minus the number of
letters $k$ is $\mu_{\overline{k}}$ if $k\leq n,$ $0$ otherwise with a
structure of cyclage graph. This structure is more complex than for type
$A_{n}.\;$In particular $\mathbf{ST}(\mu)$ decomposes into connected
components.\ These components can be isomorphic and do not necessarily contain
a row tableau. We define on $\mathbf{ST}(\mu)$ a charge statistic
$\mathrm{ch}_{n}$. Many computations allows us to conjecture that an analogue
to (\ref{th_ls}) exists for type $C_{n}$ with the charge $\mathrm{ch}_{n}%
$.\ However it seems to be impossible to derive it from our Morris type
recurrence formula.

In Section $1$ we recall the Background on Kostka-Foulkes polynomials and
crystal basis theory that we need in the sequel.\ We also summarize the basic
properties of the insertion algorithm introduced in \cite{Lec}.\ Section $2$
is devoted to the Morris type recurrence formula for type $C_{n}$ and its
applications. In Section $3$ we define the cyclage graph structure on
$\mathbf{ST}(\mu)$ and give some of its properties. Finally we introduced
$\mathrm{ch}_{n}$ in Section $4$ and conjecture that it permits to compute the
Kostka-Foulkes polynomials for type $C_{n}.$

\section{Background\label{sec1}}

\subsection{Kostka-Foulkes polynomials for type $C_{n}$}

We choose to label the Dynkin diagram of $sp_{2n}$ by%
\begin{equation}
\overset{0}{\circ}\Longrightarrow\overset{1}{\circ}-\overset{2}{\circ
}-\overset{3}{\circ}-\overset{4}{\circ}-\cdot\cdot\cdot\overset{n-1}{\circ}.
\label{DD}%
\end{equation}
The weight lattice $P_{n}$ of $C_{n}$ can be identified with $\mathbb{Z}^{n}$
equipped with the orthonormal basis $\varepsilon_{\overline{i}},$
$i=1,...,n$.\ We take for the simple roots%
\[
\alpha_{0}=2\varepsilon_{\overline{1}}\text{ and }\alpha_{i}=\varepsilon
_{\overline{i+1}}-\varepsilon_{\overline{i}}\text{, }i=1,...,n-1.
\]
Then the set of positive roots of $sp_{2n}$ is%
\[
R_{n}^{+}=\{\varepsilon_{\overline{i}}-\varepsilon_{\overline{j}}%
,\varepsilon_{\overline{i}}+\varepsilon_{\overline{j}}\text{ with }1\leq
j<i\leq n\}\cup\{2\varepsilon\overline{_{i}}\text{ with }1\leq i\leq n\}.
\]
Denote by $P_{n}^{+}$ the set of dominant weights of $sp_{2n}.$ Write
$\Lambda_{0},...,\Lambda_{n-1}$ for the fundamentals weights. Then we have
$\Lambda_{i}=\varepsilon_{\overline{n}}+\cdot\cdot\cdot+\varepsilon
_{\overline{i+1}}$, $0\leq i\leq n-1.$ Consider $\lambda\in P_{n}^{+}$ and set
$\lambda=\sum_{i=0}^{n-1}\widehat{\lambda}_{i}\Lambda_{i}$ with $\widehat
{\lambda}_{i}\in\mathbb{N}.\;$The dominant weight $\lambda$ is characterized
by the partition $(\lambda_{\overline{n}},...,\lambda_{\overline{1}})$ where
$\lambda_{\overline{i}}=\widehat{\lambda}_{0}+\cdot\cdot\cdot+\widehat
{\lambda}_{i-1},$ $i=1,...,n.$ In the sequel we will identify $\lambda$ and
$(\lambda_{\overline{n}},...,\lambda_{\overline{1}})$ by setting
$\lambda=(\lambda_{\overline{n}},...,\lambda_{\overline{1}}).$ Then
$\lambda=\lambda_{\overline{1}}\varepsilon_{\overline{1}}+\cdot\cdot
\cdot+\lambda_{\overline{n}}\varepsilon_{\overline{n}}$ that is, the
$\lambda_{i}$'s are the coordinates of $\lambda$ on the basis $(\varepsilon
_{\overline{n}},...,\varepsilon_{\overline{1}})$.\ Let $\rho$ be the half sum
of positive roots.\ We have $\rho=(n,n-1,...,1).\;$For any $\lambda\in
P_{n}^{+},$ set $\left|  \lambda\right|  =\lambda_{\overline{n}}+\cdot
\cdot\cdot+\lambda_{\overline{1}}.$

\noindent The Weyl group $W_{n}$ of $sp_{2n}$ can be regarded as the sub group
of the permutation group of $\mathcal{C}_{n}=\{\overline{n},...,\overline
{2},\overline{1},1,2,...,n\}$\ generated by $s_{i}=(i,i+1)(\overline
{i},\overline{i+1}),$ $i=1,...,n-1$ and $s_{0}=(1,\overline{1})$ where for
$a,b\in\mathcal{C}_{n}$ with $a\neq b,$ $(a,b)$ is the simple transposition
which switches $a$ and $b.$ Note that any $\sigma\in W_{n}$ verifies
$\sigma(\overline{i})=\overline{\sigma(i)}$ for $i\in\{1,...,n\}.$ We denote
by $l$ the length function corresponding to the set of generators $s_{i},$ $i=0,...n-1.$

\noindent The action of $\sigma\in W_{n}$ on $\beta=(\beta_{\overline{n}%
},...,\beta_{\overline{1}})\in P_{n}$ is given by%
\[
\sigma\cdot(\beta_{\overline{n}},...,\beta_{\overline{1}})=(\beta
_{\overline{n}}^{\sigma},...,\beta_{\overline{1}}^{\sigma})
\]
where $\beta_{\overline{i}}^{\sigma}=\beta_{\sigma(\overline{i})}$ if
$\sigma(\overline{i})\in\{\overline{1},...,\overline{n}\}$ and $\beta
_{\overline{i}}^{\sigma}=-\beta_{\sigma(i)}$ otherwise.

\noindent Let $Q_{n}^{+}$ be the set of nonnegative integral linear
combinations of positives roots.\ For any $\beta=(\beta_{\overline{n}%
},...,\beta_{\overline{1}})\in P_{n}$ we set $x^{\beta}=x_{n}^{\beta
_{\overline{n}}}\cdot\cdot\cdot x_{1}^{\beta_{\overline{1}}}$ where
$x_{1},...,x_{n}$ are fixed indeterminates$.$ The $q$-analogue of the Kostant
function partition $\mathcal{P}_{q}$ is defined by%
\[
\prod_{\alpha\in P_{n}^{+}}\dfrac{1}{1-qx^{\alpha}}=\sum_{\beta\in Q_{n}^{+}%
}\mathcal{P}_{q}(\beta)x^{\beta}\text{ and }\mathcal{P}_{q}(\beta)=0\text{ if
}\beta\notin Q_{n}^{+}.
\]

\begin{definition}
Let $\lambda,\mu\in P_{n}^{+}.$ The Kostka-Foulkes polynomial $K_{\lambda,\mu
}(q)$ is defined by
\[
K_{\lambda,\mu}(q)=\sum_{\sigma\in W_{n}}(-1)^{l(\sigma)}\mathcal{P}%
_{q}(\sigma(\lambda+\rho)-(\mu+\rho)).
\]
\end{definition}

\begin{lemma}
\label{lem_lamb1=mu1}Consider $\lambda,\mu\in P_{n}^{+}$ such that
$\lambda_{\overline{n}}=\mu_{\overline{n}}$ and write $\lambda^{\prime
}=(\lambda_{\overline{n-1}},...,\lambda_{\overline{1}}),$ $\mu^{\prime}%
=(\mu_{\overline{n-1}},...,\mu_{\overline{1}}).$ Then $K_{\lambda,\mu
}(q)=K_{\lambda^{\prime},\mu^{\prime}}(q).$
\end{lemma}

\begin{proof}
Note first that $K_{\lambda^{\prime},\mu^{\prime}}(q)$ is a Kostka-Foulkes
polynomial for type $C_{n-1}.$ We identify $W_{n-1}$ with the sub-group of
$W_{n}$ generated by $s_{i},$ $i=1,...,n-2$ and $s_{0}.$ Consider $\sigma\in
W_{n}$ such that $\sigma\notin W_{n-1}$ and set $\beta=\sigma(\lambda
+\rho)-(\mu+\rho)\in P_{n}.$ We must have $\sigma(\overline{n})\neq
\overline{n}$ since $\sigma\notin W_{n-1}.$ This implies that $\beta
=(\beta_{\overline{n}},...,\beta_{\overline{1}})$ with $\beta_{\overline{n}%
}<0$ because $\lambda_{\overline{n}}=\mu_{\overline{n}}.$ Thus $\beta\notin
Q_{n}^{+}$ and $\mathcal{P}_{q}(\beta)=0.$ This means that
\[
K_{\lambda,\mu}(q)=\sum_{\sigma\in W_{n-1}}(-1)^{l(\sigma)}\mathcal{P}%
_{q}(\sigma(\lambda^{\prime}+\rho^{\prime})-(\mu^{\prime}+\rho^{\prime})).
\]
with $\rho^{\prime}=(n-1,...,1)$ the half sum of the positive roots of the
roots system $C_{n-1}.$ Hence $K_{\lambda,\mu}(q)=K_{\lambda^{\prime}%
,\mu^{\prime}}(q).$
\end{proof}

Let $\beta\in P_{n}.\;$We set
\[
a_{\beta}=\sum_{\sigma\in W_{n}}(-1)^{l(\sigma)}(\sigma\cdot x^{\beta})
\]
where $\sigma\cdot x^{\mu}=x^{\sigma(\mu)}.$ The Schur function $s_{\beta}$ is
defined by
\[
s_{\beta}=\dfrac{a_{\beta+\rho}}{a_{\rho}}.
\]
When $\lambda\in P_{n}^{+},$ $s_{\lambda}$ is the Weyl character of
$V(\lambda)$ the finite dimensional irreducible $U_{q}(sp_{2n})$-module with
highest weight $\lambda.$ For any $\sigma\in W_{n},$ the dot action of
$\sigma$ on $\beta\in P_{n}$ is defined by $\sigma\circ\beta=\sigma\cdot
(\beta+\rho)-\rho.$ We have the following straightening law for the Schur
functions. For any $\beta\in P_{n}$ there exists a unique $\lambda\in
P_{n}^{+}$ such that $s_{\beta}=(-1)^{l(\sigma)}s_{\lambda}$ with $\sigma\in
W_{n}$ and $\lambda=\sigma\circ\beta.$ Set $\mathbb{K}=\mathbb{Z}[q,q^{-1}]$
and write $\mathbb{K}[P_{n}]$ for the $\mathbb{K}$ module generated by the
$x^{\beta}$, $\beta\in P_{n}.$ Set $\mathbb{K}[P_{n}]^{W_{n}}=\{f\in
\mathbb{K}[P_{n}],$ $\sigma\cdot f=f$ for any $\sigma\in W_{n}\}.$ Then
$\{s_{\lambda}\}$ is a basis of $\mathbb{K}[P_{n}]^{W_{n}}.$

\noindent To each positive root $\alpha,$ we associate the raising operator
$R_{\alpha}:P_{n}\rightarrow P_{n}$ defined by%
\[
R_{\alpha}(\beta)=\alpha+\beta.
\]
Given $\alpha_{1},...,\alpha_{p}$ positives roots and $\beta\in P_{n},$ we set
$(R_{\alpha_{1}}\cdot\cdot\cdot R_{\alpha_{p}})s_{\beta}=s_{R_{\alpha_{1}%
}\cdot\cdot\cdot R_{\alpha_{p}}(\beta)}.$ Composing the action of raising
operators on Schur function should be avoided in general.\ For example
$(R_{\alpha_{1}}R_{\alpha_{2}})(s_{\beta})$ is not necessarily equal to
$(R_{\alpha_{1}})(R_{\alpha_{2}}s_{\beta})$ (see example p 360 in \cite{NR}).
For all $\beta\in P_{n},$ we define the Hall-Littelwood polynomial $Q_{\beta}$
by%
\[
Q_{\beta}=\left(  \prod_{\alpha\in R_{n}^{+}}\dfrac{1}{1-qR_{\alpha}}\right)
s_{\beta}%
\]
where $\dfrac{1}{1-qR_{\alpha}}=\sum_{k=0}^{+\infty}q^{k}R_{\alpha}^{k}.$

\begin{theorem}
\label{Th_hall_kostka}\cite{mac}For any $\lambda,\mu\in P_{n}^{+},$
$K_{\lambda,\mu}(q)$ is the coefficient of $s_{\lambda}$ in $Q_{\mu}$ that
is,
\[
Q_{\mu}=\sum_{\lambda\in P_{n}^{+}}K_{\lambda,\mu}(q)s_{\lambda}.
\]
\end{theorem}

When $n=1,.$the root system $C_{1}$ can be regarded as the root system $A_{1}$
and the Kostka-Foulkes polynomial $K_{\lambda,\mu}(q)$ where $\lambda$ and
$\mu$ are partitions of length $1$ satisfies%
\begin{equation}
K_{\lambda,\mu}(q)=q^{(\left|  \lambda\right|  -\left|  \mu\right|
)/2}.\label{Kostka_n=1}%
\end{equation}

\subsection{Crystal graphs for type $C_{n}$}

Recall that crystal graphs for the $U_{q}(sp_{2n})$-modules are oriented
colored graphs with colors $i\in\{0,...,n-1\}$. An arrow $a\overset
{i}{\rightarrow}b$ means that $\widetilde{f}_{i}(a)=b$ and $\widetilde{e}%
_{i}(b)=a$ where $\widetilde{e}_{i}$ and $\widetilde{f}_{i}$ are the crystal
graph operators (for a review of crystal bases and crystal graphs see
\cite{Ka2}). A vertex $v^{0}\in B$ satisfying $\widetilde{e}_{i}(v^{0})=0$ for
any $i\in\{0,...,n-1\}$ is called a highest weight vertex. The decomposition
of $V$ into its irreducible components is reflected into the decomposition of
$B$ into its connected components. Each connected component of $B$ contains a
unique highest weight vertex.\ The crystals graphs of two isomorphic
irreducible components are isomorphic as oriented colored graphs. The action
of $\widetilde{e}_{i}$ and $\widetilde{f}_{i}$ on $B\otimes B^{\prime
}=\{b\otimes b^{\prime};$ $b\in B,b^{\prime}\in B^{\prime}\}$ is given by:%

\begin{align}
\widetilde{f_{i}}(u\otimes v)  &  =\left\{
\begin{tabular}
[c]{c}%
$\widetilde{f}_{i}(u)\otimes v$ if $\varphi_{i}(u)>\varepsilon_{i}(v)$\\
$u\otimes\widetilde{f}_{i}(v)$ if $\varphi_{i}(u)\leq\varepsilon_{i}(v)$%
\end{tabular}
\right. \label{TENS1}\\
&  \text{and}\nonumber\\
\widetilde{e_{i}}(u\otimes v)  &  =\left\{
\begin{tabular}
[c]{c}%
$u\otimes\widetilde{e_{i}}(v)$ if $\varphi_{i}(u)<\varepsilon_{i}(v)$\\
$\widetilde{e_{i}}(u)\otimes v$ if $\varphi_{i}(u)\geq\varepsilon_{i}(v)$%
\end{tabular}
\right.  \label{TENS2}%
\end{align}
where $\varepsilon_{i}(u)=\max\{k;\widetilde{e}_{i}^{k}(u)\neq0\}$ and
$\varphi_{i}(u)=\max\{k;\widetilde{f}_{i}^{k}(u)\neq0\}$. The weight of the
vertex $u$ is defined by $\mathrm{wt}(u)=\underset{i=0}{\overset{n-1}{\sum}%
}(\varphi_{i}(u)-\varepsilon_{i}(u))\Lambda_{i}$.

\noindent The following lemma is a straightforward consequence of
(\ref{TENS1}) and (\ref{TENS2}).

\begin{lemma}
\label{lem_plu_hp}Let $u\otimes v$ $\in$ $B\otimes B^{\prime}$ $u\otimes v$ is
a highest weight vertex of $B\otimes B^{\prime}$ if and only if for any
$i\in\{0,...,n-1\}$ $\widetilde{e}_{i}(u)=0$ (i.e. $u$ is of highest weight)
and $\varepsilon_{i}(v)\leq\varphi_{i}(u).$
\end{lemma}

\noindent The Weyl group $W_{n}$ acts on $B$ by:
\begin{align}
s_{i}(u)  &  =(\widetilde{f_{i}})^{\varphi_{i}(u)-\varepsilon_{i}(u)}(u)\text{
if }\varphi_{i}(u)-\varepsilon_{i}(u)\geq0,\label{actionW}\\
s_{i}(u)  &  =(\widetilde{e_{i}})^{\varepsilon_{i}(u)-\varphi_{i}(u)}(u)\text{
if }\varphi_{i}(u)-\varepsilon_{i}(u)<0.\nonumber
\end{align}
We have the equality $\mathrm{wt}(\sigma(u))=\sigma(\mathrm{wt}(u))$ for any
$\sigma\in W_{n}$ and $u\in B.$ For any $\lambda\in P_{n}^{+},$ we denote by
$B(\lambda)$ the crystal graph of $V(\lambda).$

\noindent According to (\ref{DD}) we have
\[
B(\Lambda_{n-1}):\overline{n}\overset{n-1}{\rightarrow}\overline{n-1}%
\overset{n-2}{\rightarrow}\cdot\cdot\cdot\cdot\rightarrow\overline{2}%
\overset{1}{\rightarrow}\overline{1}\overset{0}{\rightarrow}1\overset
{1}{\rightarrow}2\cdot\cdot\cdot\cdot\overset{n-2}{\rightarrow}n-1\overset
{n-1}{\rightarrow}n.
\]
Kashiwara-Nakashima's combinatorial description of the crystal graphs
$B(\lambda)$ is based on the notion of symplectic tableaux analogous for type
$C_{n}$ to semi-standard tableaux. In the sequel we use De Concini's version
of these tableaux which is equivalent to Kashiwara-Nakashima's one.

\noindent We defined a total order on $\mathcal{C}_{n}$ by setting%
\[
\mathcal{C}_{n}=\{\overline{n}<\cdot\cdot\cdot<\overline{1}<1<\cdot\cdot
\cdot<n\}.
\]
For any letter $x\in\mathcal{C}_{n}$ we set $\overline{\overline{x}}=x.$ Note
that our convention for labelling the crystal graph of the vector
representation are not those used by Kashiwara and Nakashima.\ To obtain the
original description of $B(\lambda)$ from that used in the sequel it suffices
to change each letter $k\in\{1,...,n\}$ of $\mathcal{C}_{n}$ into
$\overline{n-k+1}$ and each letter $\overline{k}\in\{\overline{1}%
,...,\overline{n}\}$ into $n-k+1.\;$The interest of this change of convention
will appear in Sections \ref{sec_cy_graph} and \ref{sec_cha}.

We identify the vertices of the crystal graph $G_{n}=\underset{l}{%
{\textstyle\bigoplus}
}B(\Lambda_{1})^{\bigotimes l}$ with the words on $\mathcal{C}_{n}$.$\;$For
any $w\in G_{n}$ we have $\mathrm{wt}(w)=d_{\overline{n}}\varepsilon
_{\overline{n}}+d_{\overline{n-1}}\varepsilon_{\overline{n-1}}\cdot\cdot
\cdot+d_{\overline{1}}\varepsilon_{\overline{1}}$ where for any $i=1,...,n$
$d_{\overline{i}}$ is the number of letters $\overline{i}$ of $w$ minus the
number of its letters $i.$ Using Formulas (\ref{TENS1}) and (\ref{TENS2}) we
obtain a simple rule to compute the action of $\widetilde{e}_{i}$,
$\widetilde{f}_{i}$ or $s_{i}$ on $w\in G_{n}$ that we will use in Section
\ref{sec_cy_graph}. Consider the subword $w_{i}$ of $w$ containing only the
letters $\overline{i+1},\overline{i},i,i+1$. Then encode in $w_{i}$ each
letter $\overline{i+1}$ or $i$ by the symbol $+$ and each letter $\overline
{i}$ or $i+1$ by the symbol $-$. Because $\widetilde{e}_{i}(+-)=\widetilde
{f}_{i}(+-)=0$ in $B(\Lambda_{n-1})\oplus B(\Lambda_{n-1})$ the factors of
type $+-$ may be ignored in $w_{i}.$ So we obtain a subword $w_{i}^{(1)}$ in
which we can ignore all the factors $+-$ to construct a new subword
$w_{i}^{(2)}$ etc... Finally we obtain a subword $\rho(w)$ of $w$ of type%
\[
\rho(w)=-^{r}+^{s}.
\]
Then we have the

\begin{claim}
\ \ \ \ \ \ \label{+-}

\begin{itemize}
\item  If $r>0,$ $\widetilde{e}_{i}(w)$ is obtained by changing the rightmost
symbol $-$ of $\rho(w)$ into its corresponding symbol $+$ (i.e. $i+1$ into $i$
and $\overline{i}$ into $\overline{i+1}$) the others letters of $w$ being
unchanged. If $r=0,$ $\widetilde{e}_{i}(w)=0.$

\item  If $s>0,$ $\widetilde{f}_{i}(w)$ is obtained by changing the leftmost
symbol $+$\ of $\rho(w)$ into their corresponding symbols $-$ (i.e. $i$ into
$i+1$ and $\overline{i+1}$ into $\overline{i}$) the others letters of $w$
being unchanged. If $s=0,\widetilde{f}_{i}(w)=0.$

\item  If $r\geq s,$ $s_{i}(w)$ is obtained by changing the $r-s$ rightmost
symbols $-$ of $\rho(w)$ into its corresponding symbol $+,$ otherwise
$s_{i}(w)$ is obtained by changing the $s-r$ leftmost symbols $+$ of $\rho(w)$
into their corresponding symbols $-.$
\end{itemize}
\end{claim}

A column on $\mathcal{C}_{n}$ is a Young diagram $C$ of column shape filled
from top to bottom by increasing letters of $\mathcal{C}_{n}$. The height
$h(C)$ of a column $C$ is the number of its letters. Set $\mathbf{C}(n,h)$ for
the set of columns of height $h$ on $\mathcal{C}_{n}$ i.e. with letters in
$\mathcal{C}_{n}$. The reading of the column $C\in\mathbf{C}(n,h)$ is the word
$\mathrm{w}((C)$ of $\mathcal{C}_{n}^{\ast}$ obtained by reading the letters
of $C$ from top to bottom.\ We will say that a column $C$ contains the pair
$(z,\overline{z})$ when $C$ contains the unbarred letter $z\geq1$ and the
barred letter $\overline{z}\leq\overline{1}$.\ Let $C_{1}$ and $C_{2}$ be two
columns. We will write $C_{1}\leq C_{2}$ when $h(C_{1})\geq h(C_{2})$ and the
rows of the tableau $C_{1}C_{2}$ weakly increase.

\begin{definition}
\label{def_admi}Let $C$ be a column on $\mathcal{C}_{n}$ and $I_{C}%
=\{z_{1}<\cdot\cdot\cdot<z_{r}\}$ the set of unbarred letters $z$ such that
the pair $(z,\overline{z})$ occurs in $C$. The column $C$ is $n$-admissible
when there exists a set of unbarred letters $J_{C}=\{t_{1}<\cdot\cdot
\cdot<t_{r}\}\subset\mathcal{C}_{n}$ such that:

\begin{itemize}
\item $t_{1}$ is the lowest letter of $\mathcal{C}_{n}$ satisfying:
$t_{1}>z_{1},t_{1}\notin C$ and $\overline{t_{1}}\notin C,$

\item  for $i=2,...,r$, $t_{i}$ is the lowest letter of $\mathcal{C}_{n}$
satisfying: $t_{i}>\max(t_{i-1,}z_{i}),$ $t_{i}\notin C$ and $\overline{t_{i}%
}\notin C.$
\end{itemize}

In this case we write:

\begin{itemize}
\item $rC$ for the column obtained from $C$ by changing $z_{i}$ into $t_{i}$
\ for each letter $z_{i}\in I_{C},$

\item $lC$ for the column obtained from $C$ by changing $\overline{z}_{i}$
into $\overline{t}_{i}$ \ for each letter $z_{i}\in I_{C}.$
\end{itemize}
\end{definition}

\noindent Consider $C=%
\begin{tabular}
[c]{|l|}\hline
$\mathtt{\bar{3}}$\\\hline
$\mathtt{\bar{2}}$\\\hline
$\mathtt{2}$\\\hline
$\mathtt{3}$\\\hline
\end{tabular}
$.\ Then $C$ is not $4$-admissible but is $5$-admissible with $rC=%
\begin{tabular}
[c]{|l|}\hline
$\mathtt{\bar{3}}$\\\hline
$\mathtt{\bar{2}}$\\\hline
$\mathtt{4}$\\\hline
$\mathtt{5}$\\\hline
\end{tabular}
$ and $lC=%
\begin{tabular}
[c]{|l|}\hline
$\mathtt{\bar{4}}$\\\hline
$\mathtt{\bar{3}}$\\\hline
$\mathtt{2}$\\\hline
$\mathtt{3}$\\\hline
\end{tabular}
.\;$As usually, we associate to each partition $\lambda=(\lambda_{\overline
{n}},...,\lambda_{\overline{1}})$ the Young diagram $Y(\lambda)$ whose $i$-th
row has length $\lambda_{\overline{n-i+1}}.\;$By definition, a $n$-symplectic
tableau $T$ of shape $\lambda$ is a filling of $Y(\lambda)$ by letters of
$\mathcal{C}_{n}$ satisfying the following conditions:

\begin{itemize}
\item  the columns $C_{i}$ of $T=C_{1}\cdot\cdot\cdot C_{s}$ are $n$-admissible,

\item  for $i=1,...,s-1:rC_{i}\leq lC_{i+1}.$
\end{itemize}

\noindent The set of $n$-symplectic tableaux will be denoted $\mathbf{ST}(n)$.
If $T=C_{1}C_{2}\cdot\cdot\cdot C_{r}\in\mathbf{ST}(n)$, the reading of $T$ is
the word $\mathrm{w}(T)=\mathrm{w}(C_{r})\cdot\cdot\cdot\mathrm{w}%
(C_{2})\mathrm{w}(C_{1})$.\ From \cite{KN} we deduce the

\begin{theorem}
\ \ \ \ \ \label{TH_KN}

\noindent\textrm{(i)}: The vertices of $B(\Lambda_{p})$ $p=0,...,n-1$ are in
one-to-one correspondence with the readings of $n$-admissible columns of
height $n-p.$

\noindent\textrm{(ii)}: The vertices of $B(\lambda)$ are in one-to-one
correspondence with the readings of the $n$-symplectic tableaux of shape
$\lambda$.
\end{theorem}

\noindent More precisely Kashiwara and Nakashima realize $B(\lambda)$ into a
tensor power $B(\Lambda_{n-1})^{\bigotimes l}$.\ Given $p=0,...,n-1,$
$B(\Lambda_{p})$ can then be identified with the connected component of
$G_{n}$ whose highest weight vertex is $b_{\overline{p}}=\overline
{n}(\overline{n-1})\cdot\cdot\cdot\overline{p+1}$.\ In this identification,
the vertices of $B(\Lambda_{p})$ are the readings of the admissible columns of
height $n-p$.\ If $\lambda=\underset{p=0}{\overset{n-1}{\sum}}\widehat
{\lambda}_{p}\Lambda_{p},$ $B(\lambda)$ is identified with the connected
component whose highest weight vertex is $b_{\lambda}=b_{\overline{n}%
}^{\otimes\widehat{\lambda}_{n}}\cdot\cdot\cdot\otimes b_{\overline{1}%
}^{\otimes\widehat{\lambda}_{1}}b_{0}^{\otimes\widehat{\lambda}_{0}}$.

\noindent By identifying $U_{q}(sp_{2(n-1)})$ with the sub-algebra of
$U_{q}(sp_{2n})$ generated by the Chevalley's generators $e_{i},f_{i}$ and
$t_{i},$ $i=0,...,n-1,$ we endow $B(\lambda)$ with a structure of crystal
graph for type $C_{n-1}.$ The decomposition of $B(\lambda)$ into its
$U_{q}(sp_{2(n-1)})$-connected components is obtained by erasing all the
arrows of color $n-1.$

\subsection{Insertion scheme for symplectic tableaux \label{par_inser}}

In \cite{Lec} we have introduced an insertion scheme for symplectic tableaux
analogous for type $C_{n}$ to the bumping algorithm on Young tableaux. Now we
are going to summarize the properties of this scheme that we shall need in
Section \ref{sec_cy_graph}.

Consider first a letter $x$ and a column $C.\;$The insertion of the letter $x$
in the $n$-admissible column $C$ is denoted $x\rightarrow C$.\ If $x$ is
strictly greater to the greatest letter of $C$ then $x\rightarrow C$ is the
column obtained by adding a box containing $x$ on bottom of $C,$ that is,
$x\rightarrow C=$%
\begin{tabular}
[c]{|l|}\hline
$C$\\\hline
$x$\\\hline
\end{tabular}
$.\;$Now suppose that $x$ is less than the greatest letter of $C$.\ Then
$x\rightarrow C$ is a symplectic tableau of two columns defined recursively as follows:

\noindent if $C=%
\begin{tabular}
[c]{|l|}\hline
$a$\\\hline
\end{tabular}
$ contains only one column then $x\rightarrow%
\begin{tabular}
[c]{|l|l|}\hline
$x$ & $a$\\\hline
\end{tabular}
$

\noindent if $C=%
\begin{tabular}
[c]{|l|}\hline
$a$\\\hline
$b$\\\hline
\end{tabular}
$ contains two letters,

\begin{enumerate}
\item $x\rightarrow%
\begin{tabular}
[c]{|l|}\hline
$a$\\\hline
$b$\\\hline
\end{tabular}
=%
\begin{tabular}
[c]{|l|l}\hline
$a$ & \multicolumn{1}{|l|}{$b$}\\\hline
$x$ & \\\cline{1-1}%
\end{tabular}
$ if $a<x\leq b$ and $b\neq\overline{a},$

\item $x\rightarrow%
\begin{tabular}
[c]{|l|}\hline
$a$\\\hline
$b$\\\hline
\end{tabular}
=%
\begin{tabular}
[c]{|l|l}\hline
$x$ & \multicolumn{1}{|l|}{$a$}\\\hline
$b$ & \\\cline{1-1}%
\end{tabular}
$ if $x\leq a<b$ and $b\neq\overline{x},$

\item $x\rightarrow%
\begin{tabular}
[c]{|c|}\hline
$\overset{\text{ \ \ \ }}{\overline{b}}$\\\hline
$b$\\\hline
\end{tabular}
=%
\begin{tabular}
[c]{|c|c}\hline
$\overset{\text{ \ \ \ }}{\overline{b+1}}$ & \multicolumn{1}{|c|}{$b+1$%
}\\\hline
$x$ & \\\cline{1-1}%
\end{tabular}
$ if $a=\overline{b}$ and $\overline{b}\leq x\leq b,$

\item $\overline{b}\rightarrow%
\begin{tabular}
[c]{|c|}\hline
$a$\\\hline
$b$\\\hline
\end{tabular}
=%
\begin{tabular}
[c]{|c|c}\hline
$\overset{\text{ \ \ \ }}{\overline{b-1}}$ & \multicolumn{1}{|c|}{$a$}\\\hline
$b-1$ & \\\cline{1-1}%
\end{tabular}
$ if $x=\overline{b}$ and $\overline{b}<a<b$.
\end{enumerate}

\noindent Consider a $n$-admissible column $C$ of height $k\geq3$ and suppose
we have defined our insertion for the $n$-admissible columns of height
$<k$.\ Set $\mathrm{w}(C)=a_{1}\cdot\cdot\cdot a_{k-1}a_{k}$ and $x\rightarrow%
\begin{tabular}
[c]{|l|}\hline
$a_{k-1}$\\\hline
$a_{k}$\\\hline
\end{tabular}
=%
\begin{tabular}
[c]{|l|l}\hline
$\delta_{k-1}$ & \multicolumn{1}{|l|}{$y$}\\\hline
$d_{k}$ & \\\cline{1-1}%
\end{tabular}
$.\ Then we have $y>a_{k-2}$ and the column $C^{\prime}$ of reading
$a_{1}\cdot\cdot\cdot a_{k-2}y$ is $n$-admissible. Write $\delta
_{k-1}\rightarrow C^{\prime}=%
\begin{tabular}
[c]{|c|c}\hline
$d_{1}$ & \multicolumn{1}{|c|}{$\ \ z$ \ \ }\\\hline
$\cdot$ & \\\cline{1-1}%
$\cdot$ & \\\cline{1-1}%
$d_{k-1}$ & \\\cline{1-1}%
\end{tabular}
$.\ We set $x\rightarrow C=$%
\begin{tabular}
[c]{|c|c}\hline
$d_{1}$ & \multicolumn{1}{|c|}{$\ \ z$ \ \ }\\\hline
$\cdot$ & \\\cline{1-1}%
$\cdot$ & \\\cline{1-1}%
$d_{k-1}$ & \\\cline{1-1}%
$d_{k}$ & \\\cline{1-1}%
\end{tabular}
. This can be pictured by%
\[
x\rightarrow%
\begin{tabular}
[c]{|c|}\hline
$a_{1}$\\\hline
$\cdot$\\\hline
$a_{k-2}$\\\hline
$a_{k-1}$\\\hline
$a_{k}$\\\hline
\end{tabular}
=%
\begin{tabular}
[c]{c|c|}\cline{2-2}%
\ \ \ \ \  & $a_{1}$\\\cline{2-2}%
& $\cdot$\\\cline{2-2}%
& $a_{k-1}$\\\cline{2-2}%
& $a_{k-1}$\\\hline
\multicolumn{1}{|c|}{$x$} & $a_{k}$\\\hline
\end{tabular}
=%
\begin{tabular}
[c]{c|c}\cline{2-2}%
& \multicolumn{1}{|c|}{$a_{1}$}\\\cline{2-2}%
& \multicolumn{1}{|c|}{$\cdot$}\\\cline{2-2}%
& \multicolumn{1}{|c|}{$a_{k-2}$}\\\hline
\multicolumn{1}{|c|}{$\delta_{k-1}$} & \multicolumn{1}{|c|}{$y$}\\\hline
\multicolumn{1}{|c|}{$d_{k}$} & \\\cline{1-1}%
\end{tabular}
=\cdot\cdot\cdot=%
\begin{tabular}
[c]{|c|c}\hline
$d_{1}$ & \multicolumn{1}{|c|}{$\ \ z$ \ \ }\\\hline
$\cdot$ & \\\cline{1-1}%
$\cdot$ & \\\cline{1-1}%
$d_{k-1}$ & \\\cline{1-1}%
$d_{k}$ & \\\cline{1-1}%
\end{tabular}
.
\]
During each step we apply one of the transformations $1$ to $4$ below. We have
proved in \cite{Lec} that $x\rightarrow C$ is then a $n$-symplectic tableau
with two columns respectively of height $h(C)$ and $1$.

\begin{example}
\label{ex_inser_col}Suppose $n=5.$

\begin{itemize}
\item $5\rightarrow%
\begin{tabular}
[c]{|l|}\hline
$\mathtt{\bar{4}}$\\\hline
$\mathtt{\bar{2}}$\\\hline
$\mathtt{2}$\\\hline
$\mathtt{3}$\\\hline
$\mathtt{4}$\\\hline
\end{tabular}
=%
\begin{tabular}
[c]{|l|}\hline
$\mathtt{\bar{4}}$\\\hline
$\mathtt{\bar{2}}$\\\hline
$\mathtt{2}$\\\hline
$\mathtt{3}$\\\hline
$\mathtt{4}$\\\hline
$\mathtt{5}$\\\hline
\end{tabular}
$.

\item $\bar{4}\rightarrow%
\begin{tabular}
[c]{|l|}\hline
$\mathtt{\bar{4}}$\\\hline
$\mathtt{\bar{2}}$\\\hline
$\mathtt{2}$\\\hline
$\mathtt{3}$\\\hline
$\mathtt{4}$\\\hline
\end{tabular}
=%
\begin{tabular}
[c]{c|c|}\cline{2-2}%
& $\mathtt{\bar{4}}$\\\cline{2-2}%
& $\mathtt{\bar{2}}$\\\cline{2-2}%
& $\mathtt{2}$\\\cline{2-2}%
& $\mathtt{3}$\\\hline
\multicolumn{1}{|c|}{$\mathtt{\bar{4}}$} & $\mathtt{4}$\\\hline
\end{tabular}
=%
\begin{tabular}
[c]{l|l}\cline{2-2}%
& \multicolumn{1}{|l|}{$\mathtt{\bar{4}}$}\\\cline{2-2}%
& \multicolumn{1}{|l|}{$\mathtt{\bar{2}}$}\\\cline{2-2}%
& \multicolumn{1}{|l|}{$\mathtt{2}$}\\\hline
\multicolumn{1}{|l|}{$\mathtt{\bar{3}}$} & \multicolumn{1}{|l|}{$\mathtt{3}$%
}\\\hline
\multicolumn{1}{|l|}{$\mathtt{3}$} & \\\cline{1-1}%
\end{tabular}
=%
\begin{tabular}
[c]{l|l}\cline{2-2}%
& \multicolumn{1}{|l|}{$\mathtt{\bar{4}}$}\\\cline{2-2}%
& \multicolumn{1}{|l|}{$\mathtt{\bar{2}}$}\\\hline
\multicolumn{1}{|l|}{$\mathtt{\bar{2}}$} & \multicolumn{1}{|l|}{$\mathtt{2}$%
}\\\hline
\multicolumn{1}{|l|}{$\mathtt{2}$} & \\\cline{1-1}%
\multicolumn{1}{|l|}{$\mathtt{3}$} & \\\cline{1-1}%
\end{tabular}
=%
\begin{tabular}
[c]{l|l}\cline{2-2}%
& \multicolumn{1}{|l|}{$\mathtt{\bar{4}}$}\\\hline
\multicolumn{1}{|l|}{$\mathtt{\bar{3}}$} & \multicolumn{1}{|l|}{$\mathtt{3}$%
}\\\hline
\multicolumn{1}{|l|}{$\mathtt{\bar{2}}$} & \\\cline{1-1}%
\multicolumn{1}{|l|}{$\mathtt{2}$} & \\\cline{1-1}%
\multicolumn{1}{|l|}{$\mathtt{3}$} & \\\cline{1-1}%
\end{tabular}
=%
\begin{tabular}
[c]{|l|l}\hline
$\mathtt{\bar{4}}$ & \multicolumn{1}{|l|}{$\mathtt{3}$}\\\hline
$\mathtt{\bar{3}}$ & \\\cline{1-1}%
$\mathtt{\bar{2}}$ & \\\cline{1-1}%
$\mathtt{2}$ & \\\cline{1-1}%
$\mathtt{3}$ & \\\cline{1-1}%
\end{tabular}
.$
\end{itemize}
\end{example}

\noindent\textbf{Remarks:}

\noindent$\mathrm{(i)}\mathbf{:}$ If $C$ is $n$-admissible and $x\rightarrow
C$ is a column, then this column is $(n+1)$-admissible but not necessarily $n$-admissible.

\noindent$\mathrm{(ii)}\mathbf{:}$ From transformations $1$ to $4$ we obtain
the plactic relations%
\begin{equation}
abx=\left\{
\begin{tabular}
[c]{l}%
$bax$ if $a<x\leq b$ and $b\neq\overline{a}$\\
$axb$ if $x\leq a<b$ and $b\neq\overline{x}$%
\end{tabular}
\right.  \text{ and }\left\{
\begin{tabular}
[c]{l}%
$\overline{b}bx=(b+1)(\overline{b+1})x$ if $\overline{b}\leq x\leq b$\\
$ab\overline{b}=a(\overline{b-1})(b-1)$ if $\overline{b}<a<b$%
\end{tabular}
\right.  \label{placrela}%
\end{equation}
introduced in \cite{Lec}. These relations are not sufficient to define a
plactic monoid for type $C_{n}.\;$We need a contraction relation which permits
to obtain a $n$-admissible column from a non $n$-admissible one. Let
$C^{\prime}=$%
\begin{tabular}
[c]{|l|}\hline
$C$\\\hline
$x$\\\hline
\end{tabular}
be a non $n$-admissible column on $\mathcal{C}_{n}$ such that $C$ is
$n$-admissible and $x$ a letter.\ In this case we can prove that there exists
an unbarred letter $z$ maximal such that the pair $(z,\overline{z})$ occurs in
$C^{\prime}$ and
\[
\mathrm{card}\{t\in\mathcal{C}_{n},\left|  t\right|  \geq z\}>n-z+1.
\]
Write $D$ for the column obtained by erasing the pair $(z,\overline{z})$ in
$C^{\prime}.$ Then $D$ is $n$-admissible and the contraction relation is
defined by%
\begin{equation}
\mathrm{w}(C^{\prime})\equiv\mathrm{w}(D)\text{.} \label{relcontra}%
\end{equation}
In fact this last relation is not needed to define the cocyclage in Section
\ref{sec_cy_graph}. We will denote $\equiv_{n}$ the congruence obtained by
identifying words of $\mathcal{C}_{n}^{\ast}$ which are equal up to relations
(\ref{placrela}).

\noindent$\mathrm{(iii)}\mathbf{:}$ When $x\rightarrow C=%
\begin{tabular}
[c]{|l|l|}\hline
$C^{\prime}$ & $y$\\\hline
\end{tabular}
$ is a tableau of two columns we have%
\begin{equation}
h(C)=h(C^{\prime})\text{ and }C^{\prime}\leq C. \label{C'infC}%
\end{equation}

Now we can define the insertion $x\rightarrow T$ of the letter $x$ in the
$n$-symplectic tableau $T$.\ Write $T=C_{1}\cdot\cdot\cdot C_{r}$ where
$C_{i},$ $i=1,...,r$ are the $n$-admissible columns of $T.\;$If $x$ is
strictly greater to the greatest letter of $C_{1}$ then $x\rightarrow T$ is
the tableau obtained by adding a box containing $x$ on bottom of $C_{1}%
$.\ Then $x\rightarrow T$ belongs to $\mathbf{ST}(n+1)\;$but not to
$\mathbf{ST}(n)$ in general since its first column may be non $n$-admissible.
Otherwise write $x\rightarrow C=%
\begin{tabular}
[c]{|l|l|}\hline
$C_{1}^{\prime}$ & $y$\\\hline
\end{tabular}
$ where $C_{1}^{\prime}$ is an admissible column of height $h(C_{1})$ and $y$
a letter.\ Then $x\rightarrow T=C_{1}^{\prime}(y\rightarrow C_{2}\cdot
\cdot\cdot C_{r})$ that is, $x\rightarrow T$ is the juxtaposition of
$C_{1}^{\prime}$ with the tableau obtained by inserting $y$ in the tableau
$C_{2}\cdot\cdot\cdot C_{r}.$ In this case $x\rightarrow T$ is a
$n$-symplectic tableau \cite{Lec}.

\begin{example}
Suppose $n=3.\;$Then $2\rightarrow\left(  1\rightarrow%
\begin{tabular}
[c]{|l|ll}\hline
$\mathtt{\bar{1}}$ & $\mathtt{1}$ & \multicolumn{1}{|l|}{$\mathtt{2}$}\\\hline
$\mathtt{1}$ & $\mathtt{2}$ & \multicolumn{1}{|l}{}\\\cline{1-2}%
$\mathtt{3}$ &  & \\\cline{1-1}%
\end{tabular}
\right)  =%
\begin{tabular}
[c]{|l|l|ll}\hline
$\mathtt{\bar{2}}$ & $\mathtt{1}$ & $\mathtt{2}$ &
\multicolumn{1}{|l|}{$\mathtt{2}$}\\\hline
$\mathtt{1}$ & $\mathtt{2}$ &  & \\\cline{1-2}\cline{2-2}%
$\mathtt{2}$ & $\mathtt{3}$ &  & \\\cline{1-2}%
\end{tabular}
$.
\end{example}

\noindent\textbf{Remarks: }

\noindent$\mathrm{(i):}$ The insertion scheme described below do not suffice
to define a complete insertion algorithm for the $n$-symplectic tableaux since
the first column $C^{\prime}$ of $x\rightarrow T$ may be not $n$-admissible
when $x$ is greater than the greatest letter of $C_{1}.\;$To obtain a complete
insertion algorithm we have to apply relation (\ref{relcontra}) to $C^{\prime
}$.\ This give a column $D$ of reading $x_{1}\cdot\cdot\cdot x_{p}.\;$Finally
we compute successively the insertions $x_{p}(\rightarrow x_{p-1}\cdot
\cdot\cdot(x_{1}\rightarrow C_{2}\cdot\cdot\cdot C_{r}))$. In the sequel we
only use insertion algorithm without the contraction relation (\ref{relcontra}).

\noindent$\mathrm{(ii):}$ To each $w=x_{1}\cdot\cdot\cdot x_{r}\in
\mathcal{C}_{n}^{\ast}$ of length $r$ we can associate recursively a
symplectic tableau $P(w)$ by setting $P(w)=$%
\begin{tabular}
[c]{|l|}\hline
$x_{1}$\\\hline
\end{tabular}
if $r=1$ and $P(w)=x_{r}\rightarrow P(x_{1}\cdot\cdot\cdot x_{r-1})$
otherwise. If $P(w)$ belongs to $\mathbf{ST}(m)$ with $m\geq n,$ results of
\cite{Lec} implies that $P(w)\equiv_{m}w.$ Moreover $P(w)$ is the unique
$m$-symplectic tableau with this property.\ Denote by $\sim_{m}$ the
equivalence relation defining on the vertices of $G_{m}$ by $w_{1}\sim
_{m}w_{2}$ if and only if $w_{2}$ and $w_{2}$ belong to the same connected
component of $G_{m}.\;$Given two words $w_{1}$ and $w_{2}$ such that
$P(w_{1})$ and $P(w_{2})$ belong $\mathbf{ST}(m)$ we have the equivalences
\begin{equation}
w_{1}\equiv_{m}w_{2}\Longleftrightarrow P(w_{1})=P(w_{2})\Longleftrightarrow
P(w_{1})\sim_{m}P(w_{2}). \label{char_placti}%
\end{equation}
Moreover we have for any $\sigma\in W_{m}$%
\begin{equation}
P(\sigma(w))=\sigma(P(w)). \label{rq_p_com_w}%
\end{equation}

\noindent$\mathrm{(iii):}$ The insertion algorithm is reversible in the sense
that if we know the tableau $T^{\prime}$ such that $x\rightarrow T=T^{\prime}$
and the shape of $T$ we can recover the tableau $T$ and the letter $x$.\ This
follows from the fact that the transformations $1$ to $4$ are reversible. More
precisely, $T^{\prime}$ has one box more than $T.\;$Let $y$ be the letter
belonging to that box. Then if we apply transformations $1$ to $4$ from right
to left starting from $y,$ we recover $T$ and $x.$

\noindent$\mathrm{(iv):}$ In Section \ref{sec_cy_graph}, we will need to find
for a fixed tableau $T^{\prime}$ all the pairs $(x,T)$ where $T$ is a
symplectic tableau and $x$ a letter such that $x\rightarrow T=T^{\prime}$. The
outside corners of the tableau $T^{\prime}$ are the boxes $c$ of $T^{\prime}$
such that there is no box down and to the right of $c$ in $T^{\prime}.\;$By
$\mathrm{(iii)}$ the pairs $(x,T)$ are obtained by applying the reverse
insertion algorithm to the outside corners of $T^{\prime}.$

\begin{example}
Suppose $n=3$ and $T^{\prime}=%
\begin{tabular}
[c]{|l|ll}\hline
$\mathtt{\bar{2}}$ & $\mathtt{1}$ & \multicolumn{1}{|l|}{$\mathbf{2}$}\\\hline
$\mathtt{1}$ & $\mathbf{2}$ & \multicolumn{1}{|l}{}\\\cline{1-2}%
$\mathbf{3}$ &  & \\\cline{1-1}%
\end{tabular}
.\;$Then by applying reverse insertion algorithm to each outside corners of
$T^{\prime}$ we obtain the pairs $\left(  3,%
\begin{tabular}
[c]{|l|l|l}\hline
$\mathtt{\bar{2}}$ & $\mathtt{1}$ & \multicolumn{1}{|l|}{$\mathtt{2}$}\\\hline
$\mathtt{1}$ & $\mathtt{2}$ & \\\cline{1-2}%
\end{tabular}
\right)  ,$ $\left(  1,%
\begin{tabular}
[c]{|l|ll}\hline
$\mathtt{\bar{1}}$ & $\mathtt{1}$ & \multicolumn{1}{|l|}{$\mathtt{2}$}\\\hline
$\mathtt{1}$ &  & \\\cline{1-1}%
$\mathtt{3}$ &  & \\\cline{1-1}%
\end{tabular}
\right)  $ and $\left(  1,%
\begin{tabular}
[c]{|l|l}\hline
$\mathtt{\bar{1}}$ & \multicolumn{1}{|l|}{$\mathtt{1}$}\\\hline
$\mathtt{1}$ & \multicolumn{1}{|l|}{$\mathtt{2}$}\\\hline
$\mathtt{3}$ & \\\cline{1-1}%
\end{tabular}
\right)  .$
\end{example}

\section{Morris type recurrence formula}

In this section we introduce a recurrence formula for computing Kostka
polynomials analogous for type $C_{n}$ to Morris recurrence formula.\ It
allows to explain the Kostka polynomials for type $C_{n}$ as combinations of
Kostka polynomials for type $C_{n-1}.$ We embed type $C_{n-1}$ in type $C_{n}$
by identifying $U_{q}(sp_{2(n-1)})$ with the sub-algebra of $U_{q}(sp_{2n})$
generated by the Chevalley operators $e_{i},f_{i}$ and $t_{i},$ $i=0,...n-2.$
The weight lattice $P_{n-1}$ of $U_{q}(sp_{2n})$ is the $\mathbb{Z}$-lattice
generated by the $\varepsilon_{\overline{i}},$ $i=1,...,n-1$ and $P_{n-1}%
^{+}=P_{n}^{+}\cap P_{n-1}$ is the set of dominant weights. The Weyl group
$W_{n-1}$ is the sub-group of $W_{n}$ generated by the $s_{i},$ $i=0,...n-2$
and we have $R_{n-1}^{+}=R_{n}\cap P_{n-1}.$

\noindent Given any positive integer $r,$ write $(r)_{n}$ for the row
partition $(p,0,...0)$ of length $n$.\ To obtain our recurrence formula we
need to describe the decomposition $B(\gamma)\otimes B((r)_{n})$ with
$\gamma\in P_{n}^{+}$ and $r>0$ an integer into its irreducible
components.\ This is analogous for type $C_{n}$ to Pieri rule.

\subsection{Pieri rule for type $C_{n}$}

Let $\gamma=(\gamma_{\overline{n}},...,\gamma_{\overline{1}})\in P_{n}^{+}.$
By Theorem \ref{TH_KN}, the vertices of $B((r)_{n})$ are the words
\[
L=(n)^{k_{n}}\cdot\cdot\cdot(2)^{k_{2}}(1)^{k_{1}}(\overline{1})^{k_{\bar{1}}%
}(\overline{2})^{k_{\bar{2}}}\cdot\cdot\cdot(\overline{n})^{k_{\bar{n}}}%
\]
where $k_{\overline{i}},k_{i}$ are positive integers, $(x)^{k}$ means that the
letter $x$ is repeated $k$ times in $L$ and $k_{\overline{1}}+\cdot\cdot
\cdot+k_{\overline{n}}+k_{1}+\cdot\cdot\cdot+k_{n}=r.$ Let $b_{\gamma}$ be the
highest weight vertex of $B(\gamma).$

\begin{lemma}
$b_{\gamma}\otimes L$ is a highest weight vertex of $B(\gamma)\otimes
B((r)_{n})$ if and only if the following conditions holds:

\noindent$\mathrm{(i)}:$ $\gamma_{\overline{i}}-k_{i}\geq\gamma_{\overline
{i-1}}$ for $i=2,...,n$ and $\gamma_{\overline{1}}-k_{1}\geq0,$

\noindent$\mathrm{(ii)}:$ $\gamma_{\overline{i}}-k_{i}+k_{\overline{i}}%
\leq\gamma_{\overline{i+1}}-k_{i+1}$ for $i=1,...,n-1.$
\end{lemma}

\begin{proof}
By Lemma \ref{lem_plu_hp}, $b_{\gamma}\otimes L$ is a highest weight vertex if
and only if for any $m=1,...r,$ each vertex $b_{\gamma}\otimes L_{m}$ (where
$L_{m}$ is the word obtained by reading the $m$ leftmost letters of $L$) is a
highest weight vertex. It means that $(\gamma_{\overline{n}}-k_{n}%
,..,\gamma_{\overline{s}}-k_{s},\gamma_{\overline{s-1}},...,\gamma
_{\overline{1}})$ and $(\gamma_{\overline{n}}-k_{n},..,\gamma_{\overline{t+1}%
}+k_{t+1},\gamma_{\overline{t}}-k_{t}+k_{\overline{t}},...,\gamma
_{\overline{1}}+k_{1}-k_{\overline{1}})$ are partitions respectively for
$s=n,...,1$ and $t=1,...,n-1.$ This is equivalent to the conditions%
\[
\left\{
\begin{tabular}
[c]{l}%
$\gamma_{\overline{s}}-k_{s}\geq\gamma_{\overline{s-1}}\text{ for }s=n,...,2$
and $\gamma_{\overline{1}}-k_{\overline{1}}\geq0$\\
$\gamma_{\overline{t}}-k_{t}+k_{\overline{t}}\leq\gamma_{\overline{t+1}%
}-k_{k+1}$ for $t=1,...,n-1$%
\end{tabular}
\right.  .
\]
\end{proof}

\begin{corollary}
\label{cor_pieri}$B(\gamma)\otimes B((r)_{n})=\underset{\lambda\in P_{n}^{+}%
}{\bigoplus}B(\lambda)^{\oplus n_{\lambda}}$ where $n_{\lambda}$ is the number
of vertices $L\in B((r)_{n})$ such that

\noindent$\mathrm{(i)}:$ $k_{\overline{i}}-k_{i}=\lambda_{\overline{i}}%
-\gamma_{\overline{i}}$ for $i=1,....,n$,

\noindent$\mathrm{(ii)}:$ $\lambda_{\overline{i}}\leq\lambda_{\overline{i+1}%
}-k_{\overline{i+1}}$ for $i=1,...,n-1$,

\noindent$\mathrm{(iii)}:$ $\lambda_{\overline{i}}-k_{\overline{i}}\geq
\lambda_{\overline{i-1}}+k_{i-1}-k_{\overline{i-1}}$ for $i=2,...,n$ and
$\lambda_{\overline{1}}-k_{\overline{1}}\geq0.$
\end{corollary}

\begin{proof}
The multiplicity $n_{\lambda}$ is equal to the number of highest weight
vertices $b_{\gamma}\otimes L\in B(\gamma)\otimes B((r)_{n})$ of weight
$\lambda.$ The condition $\mathrm{wt}(b_{\gamma}\otimes L)=\lambda$ is
equivalent to
\[
\gamma_{\overline{i}}-k_{i}+k\overline{_{i}}=\lambda_{\overline{i}}\text{ for
}i=1,...,n
\]
which gives $\mathrm{(i).}$ The assertions $\mathrm{(ii)}$ and $\mathrm{(iii)}%
$ are respectively obtained by replacing for any $i,$ $\gamma_{\overline{i}}$
by $\lambda_{\overline{i}}+k_{i}-k_{\overline{i}}$ in assertions
$\mathrm{(ii)}$ and $\mathrm{(i)}$ of the previous lemma.
\end{proof}

\noindent Note that $B(\gamma)\otimes B((r)_{n})$ is not multiplicity free in general.

\subsection{Recurrence formula}

\begin{theorem}
\label{Th_Morris}Let $\mu\in P_{n}^{+}$ with $n\geq2$ and write $\mu
=(\mu_{\overline{n}},\mu^{\prime})$ where $\mu_{\overline{n}}$ is the first
part of $\mu$ and $\mu^{\prime}=(\mu_{\overline{n-1}},...,\mu_{\overline{1}%
})\in P_{n-1}^{+}.$ Then
\begin{equation}
Q_{\mu}=\sum_{\gamma\in P_{n-1}^{+}}\sum_{r=0}^{+\infty}\sum_{m=0}^{+\infty
}q^{m+r}\sum_{\lambda\in B(\gamma)\otimes B((r)_{n-1})}K_{\lambda,\mu^{\prime
}}(q)s_{(\mu_{\overline{n}}+r+2m,\gamma)}\label{Morrisrec}%
\end{equation}
\end{theorem}

\begin{proof}
We start from $Q_{\mu}=\left(  \prod_{\alpha\in R_{n}^{+}}\dfrac
{1}{1-qR_{\alpha}}\right)  s_{\mu}.$ By Proposition 3.5 of \cite{NR} we can
write
\[
Q_{\mu}=\left(  \underset{\alpha\notin R_{n-1}^{+}}{\prod_{\alpha\in R_{n}%
^{+}}}\dfrac{1}{1-qR_{\alpha}}\right)  \left[  \left(  \underset{\alpha\in
R_{n-1}^{+}}{\prod}\dfrac{1}{1-qR_{\alpha}}\right)  s_{\mu}\right]  .
\]
Then by applying Theorem \ref{Th_hall_kostka}, we obtain%
\begin{equation}
Q_{\mu}=\left(  \underset{\alpha\notin R_{n-1}^{+}}{\prod_{\alpha\in R_{n}%
^{+}}}\dfrac{1}{1-qR_{\alpha}}\right)  \left(  \sum_{\lambda\in P_{n-1}^{+}%
}K_{\lambda,\mu^{\prime}}(q)s_{(\mu_{\overline{n}},\lambda)}\right)
.\label{for_pro_q_mu}%
\end{equation}
Set $R_{\overline{i}}=R_{\varepsilon_{\overline{n}}-\varepsilon_{\overline{i}%
}}$ for $i=1,...,n-1$ and $R_{i}=R_{\varepsilon_{\overline{n}}+\varepsilon
_{\overline{i}}}$ for $i=1,...,n.$ Recall that for any $\beta\in P_{n-1},$
$R_{\overline{i}}(\beta)=\beta+\varepsilon_{\overline{n}}-\varepsilon
_{\overline{i}}$ and $R_{i}(\beta)=\beta+\varepsilon_{\overline{n}%
}+\varepsilon_{\overline{i}}.$ Then (\ref{for_pro_q_mu}) implies%
\begin{multline*}
Q_{\mu}=\sum_{\lambda\in P_{n-1}^{+}}K_{\lambda,\mu^{\prime}}(q)\times\\
\left(  \sum_{r=0}^{+\infty}\sum_{m=0}^{+\infty}\sum_{k_{\overline{1}}%
+\cdot\cdot\cdot+k_{\overline{n}}+k_{1}+\cdot\cdot\cdot+k_{n}=r}q^{m+r}%
(R_{n})^{m}(R_{1})^{k_{1}}(R_{\overline{1}})^{k_{\overline{1}}}\cdot\cdot
\cdot(R_{n-1})^{k_{n-1}}(R_{\overline{n-1}})^{k_{\overline{n-1}}}%
s_{(\mu_{\overline{n}},\lambda)}\right)  .
\end{multline*}%
\[
Q_{\mu}=\sum_{r=0}^{+\infty}\sum_{m=0}^{+\infty}q^{m+r}\sum_{\lambda\in
P_{n-1}^{+}}K_{\lambda,\mu^{\prime}}(q)\sum_{k_{\overline{1}}+\cdot\cdot
\cdot+k_{\overline{n}}+k_{1}+\cdot\cdot\cdot+k_{n}=r}s_{(\mu_{\overline{n}%
}+r+2m,\lambda_{\overline{n-1}}+k_{n-1}-k_{\overline{n-1}},\cdot\cdot
\cdot,\lambda_{\overline{1}}+k_{1}-k_{\overline{1}})}.
\]
Fix $\lambda,m$ and $r$ and consider%
\[
S=\sum_{k_{\overline{1}}+\cdot\cdot\cdot+k_{\overline{n}}+k_{1}+\cdot
\cdot\cdot+k_{n}=r}s_{(\mu_{\overline{n}}+r+2m,\lambda_{\overline{n-1}%
}+k_{n-1}-k_{\overline{n-1}},\cdot\cdot\cdot,\lambda_{\overline{1}}%
+k_{1}-k_{\overline{1}})}.
\]
Set $\gamma=(\lambda_{\overline{n-1}}+k_{n-1}-k_{\overline{n-1}}%
,...,\lambda_{\overline{1}}+k_{1}-k_{\overline{1}}).$

\noindent$\mathrm{1:}$\textrm{ }Suppose first that there exits $i\in
\{1,...,n-2\}$ such that $\lambda_{\overline{i}}>\lambda_{\overline{i+1}%
}-k_{\overline{i+1}}.$ Set $\widetilde{\gamma}=s_{i}\circ\gamma$ that is
\[
\widetilde{\gamma}=s_{i}(\gamma_{\overline{n-1}}+n-1,...,\gamma_{\overline
{i+1}}+n-i+1,\gamma_{\overline{i}}+n-i,...,\gamma_{\overline{1}}%
+1)-(n-1,...,1).
\]
Then $\gamma_{\overline{s}}=\widetilde{\gamma}_{\overline{s}}$ for $s\neq
i+1,i$, $\widetilde{\gamma}_{\overline{i+1}}=\gamma_{\overline{i}}-1$ and
$\widetilde{\gamma}_{\overline{i}}=,\gamma_{\overline{i+1}}+1$ that is
\[
\left\{
\begin{tabular}
[c]{l}%
$\widetilde{\gamma}_{\overline{i+1}}=\lambda_{\overline{i}}+k_{i}%
-k_{\overline{i}}-1$\\
$\widetilde{\gamma}_{\overline{i}}=\lambda_{\overline{i+1}}+k_{i+1}%
-k_{\overline{i+1}}+1$%
\end{tabular}
\right.  .
\]
Write $\widetilde{k}_{i+1}=k_{i}$, $\widetilde{k}_{i}=k_{i+1}$, $\widetilde
{k}_{\overline{i+1}}=\lambda_{\overline{i+1}}-\lambda_{\overline{i}%
}+k_{\overline{i}}+1$ and $\widetilde{k}_{\overline{i}}=\lambda_{\overline{i}%
}-\lambda_{\overline{i+1}}+k_{\overline{i+1}}-1.$ To make our notation
homogeneous set $\widetilde{k}_{t}=k_{t}$ for any $t\neq i,i+1,\overline
{i},\overline{i+1}.$ Then $\lambda_{\overline{i}}>\lambda_{\overline{i+1}%
}-\widetilde{k}_{\overline{i+1}}.$ We have $\widetilde{k}_{\overline{i+1}}%
\geq0$ and $\widetilde{k}_{\overline{i}}=\lambda_{\overline{i}}-\lambda
_{\overline{i+1}}+k_{\overline{i+1}}-1\geq0$ since $\lambda_{\overline{i}%
}>\lambda_{\overline{i+1}}-k_{\overline{i+1}}$.$\;$Moreover $\widetilde
{k}_{\overline{1}}+\cdot\cdot\cdot+\widetilde{k}_{\overline{n}}+\widetilde
{k}_{1}+\cdot\cdot\cdot+\widetilde{k}_{n}=r$ and for any $s\in\{1,...,n-2\}$%
\[
\widetilde{\gamma}_{\overline{s}}=\lambda_{\overline{s}}+\widetilde{k}%
_{s}-\widetilde{k}_{\overline{s}}.
\]

\noindent$\mathrm{2:}$ Suppose that $\lambda_{\overline{i}}\leq\lambda
_{\overline{i+1}}-k_{\overline{i+1}}$ for any $i=1,...,n-2$ and $\lambda
_{\overline{1}}-k_{\overline{1}}<0$. Set $\widetilde{\gamma}=s_{0}\circ
\gamma.$ Then $\gamma_{\overline{s}}=\widetilde{\gamma}_{\overline{s}}$ for
$s\neq1$ and $\widetilde{\gamma}_{\overline{1}}=-\gamma_{\overline{1}}%
-k_{1}+k_{\overline{1}}-2$. Write $\widetilde{k}_{i}=k_{i},$ $\widetilde
{k}_{\overline{i}}=k_{\overline{i}}$ for any $i=2,...,n-1$ and set
$\widetilde{k}_{1}=k_{\overline{1}}-\lambda_{\overline{1}}-1,$ $\widetilde
{k}_{\overline{1}}=k_{1}+\lambda_{\overline{1}}+1.$ We have $\widetilde{k}%
_{1}\geq0$, $\lambda_{\overline{i}}\leq\lambda_{\overline{i+1}}-\widetilde
{k}_{\overline{i+1}}$ for any $i=1,...,n-2$ and $\lambda_{\overline{1}%
}-\widetilde{k}_{\overline{1}}<0.\;$Moreover $\widetilde{k}_{\overline{1}%
}+\cdot\cdot\cdot+\widetilde{k}_{\overline{n}}+\widetilde{k}_{1}+\cdot
\cdot\cdot+\widetilde{k}_{n}=r.$

\noindent$\mathrm{3:}$\textrm{ }Now suppose that $\lambda_{\overline{s}}%
\leq\lambda_{\overline{s+1}}-k_{\overline{s+1}}$ for any $s\in\{1,...,n-2\}$,
$\lambda_{\overline{1}}-k_{\overline{1}}\geq0$ and there exists $i\in
\{1,...,n-2\}$ such that $\lambda_{\overline{i+1}}-k_{\overline{i+1}}%
<\lambda_{\overline{i}}+k_{i}-k_{\overline{i}}.$ Define $\widetilde{\gamma
}=s_{i}\circ\gamma$ as above.\ Set $\widetilde{k}_{\overline{i+1}%
}=k_{\overline{i+1}}$, $\widetilde{k}_{\overline{i}}=k_{\overline{i}}$,
$\widetilde{k}_{i+1}=\lambda_{\overline{i}}-\lambda_{\overline{i+1}%
}-k_{\overline{i}}+k_{i}+k_{\overline{i+1}}-1$ and $\widetilde{k}_{i}%
=(\lambda_{\overline{i+1}}-\lambda_{\overline{i}}-k_{\overline{i+1}}%
)+k_{i+1}+k_{\overline{i}}+1.$ Write $\widetilde{k}_{t}=k_{t}$ for any $t\neq
i,i+1,\overline{i},\overline{i+1}.$ We obtain $\widetilde{k}_{i}\geq0$ and
$\widetilde{k}_{i+1}\geq0$ since $\lambda_{\overline{i}}\leq\lambda
_{\overline{i+1}}-k_{\overline{i+1}}$ and $\lambda_{\overline{i+1}%
}-k_{\overline{i+1}}<\lambda_{\overline{i}}+k_{i}-k_{\overline{i}}.$ Since
$\widetilde{k}_{\overline{s}}=k_{\overline{s}}$ for any $s=1,...,n-1,$ we have
$\lambda_{\overline{s}}\leq\lambda_{\overline{s+1}}-\widetilde{k}%
_{\overline{s+1}}$ for any $s\in\{1,...,n-2\}$ and $\lambda_{\overline{1}%
}-\widetilde{k}_{\overline{1}}\geq0.\;$Moreover the assertion $\lambda
_{\overline{i+1}}-\widetilde{k}_{\overline{i+1}}<\lambda_{\overline{i}%
}+\widetilde{k}_{i}-\widetilde{k}_{\overline{i}}$ holds since it is equivalent
to $0<k_{i+1}+1.$ Finally $\widetilde{k}_{\overline{1}}+\cdot\cdot
\cdot+\widetilde{k}_{\overline{n}}+\widetilde{k}_{1}+\cdot\cdot\cdot
+\widetilde{k}_{n}=r$ and for any $s\in\{1,...,n-2\}$%
\[
\widetilde{\gamma}_{\overline{s}}=\lambda_{\overline{s}}+\widetilde{k}%
_{s}-\widetilde{k}_{\overline{s}}.
\]

\noindent Denote by $E_{1}$ $E_{2}$ and $E_{3}$ the sets of multi-indices
$(k_{\overline{1}},...,k_{\overline{n}},k_{1},...,k_{n})$ such that
$k_{\overline{1}}+\cdot\cdot\cdot+k_{\overline{n}}+k_{1}+\cdot\cdot\cdot
+k_{n}=r$ and satisfying respectively the assertions $\mathrm{1}$\textrm{,}
$\mathrm{2,}$ $\mathrm{3}$\textrm{.\ }Let\textrm{ }$\chi$ be the map defined
on $E_{1}\cup E_{2}\cup E_{3}$ by
\[
\chi(\gamma)=\widetilde{\gamma}.
\]
Then by the above arguments $\chi$ is a bijection which verifies $\chi
(E_{i})=E_{i}$ for $i=1,2,3$.\ Now the pairing $\gamma\longleftrightarrow
\widetilde{\gamma}$ provides the cancellation of all the $s_{\gamma}$ with
$\gamma=(\lambda_{\overline{n-1}}+k_{n-1}-k_{\overline{n-1}},...,\lambda
_{\overline{1}}+k_{1}-k_{\overline{1}})$ such that $(k_{\overline{1}%
},...,k_{\overline{n}},k_{1},...,k_{n})\in E_{1}\cup E_{2}\cup E_{3}$
appearing in $S.$ Indeed $s_{(\mu_{\overline{n}}+r+2m,\gamma)}=-s_{(\mu
_{\overline{n}}+r+2m,\widetilde{\gamma})}.$ By Corollary \ref{cor_pieri} it
means that
\[
S=\underset{\lambda\in B(\gamma)\otimes B((r)_{n-1})}{\underset{\gamma\in
P_{n-1}^{+}}{\sum}}s_{(\mu_{\overline{n}}+r+2m,\gamma)}%
\]
and the theorem is proved.
\end{proof}

\noindent Note that the theorem is also true for $n=2.$ In this case $R_{n-1}$
is the set of positive roots of the root system $A_{1}.$

\begin{corollary}
Let $\nu,\mu\in P_{n}^{+}$ such that $\mu_{\overline{n}}\geq\nu_{\overline
{n-1}}.\;$Set $l=\nu_{\overline{n}}-\mu_{\overline{n}}\geq0$ and $\nu^{\prime
}=(\nu_{\overline{n-1}},...,\nu_{\overline{1}}).\;$Then
\[
K_{\nu,\mu}(q)=\sum_{r+2m=l}q^{r+m}\sum_{\lambda\in B(\nu^{\prime})\otimes
B((r)_{n-1})}K_{\lambda,\mu^{\prime}}(q).
\]
\end{corollary}

\begin{proof}
Let $m,r$ be two integers such that $(\mu_{\overline{n}}+r+2m,\nu^{\prime
})=\nu.\;$Then $l=r+2m.\;$Consider $s_{(\mu_{\overline{n}}+r_{\ast}+2m_{\ast
},\gamma)}$ and $s_{(\mu_{\overline{n}}+r+2m,\nu^{\prime})}$ appearing in
(\ref{Morrisrec}).\ Suppose that there exists $\sigma\in W_{n}$ such that
$(\mu_{\overline{n}}+r_{\ast}+2m_{\ast},\gamma)=\sigma\circ(\mu_{\overline{n}%
}+r+2m,\nu^{\prime})=\sigma(\mu_{\overline{n}}+r+2m+n,\nu_{\overline{n-1}%
}+n-1,...,\nu_{\overline{1}}+1)-(n,...,1).$ We can not have $\sigma
(p)=\overline{n}$ with $p\in\{1,...,n\}$ otherwise $\mu_{\overline{n}}%
+r_{\ast}+2m_{\ast}<0.\;$Set $\sigma(\overline{p})=\overline{n}$ with
$p\in\{1,...,n\}.$ If $p<n$ we must have%
\[
\nu_{\overline{p}}+p=\mu_{\overline{n}}+r_{\ast}+2m_{\ast}+n.
\]
Thus $\nu_{\overline{p}}=\mu_{\overline{n}}+r_{\ast}+2m_{\ast}+n-p>\mu
_{\overline{n}}$ which contradicts the hypothesis $\mu_{\overline{n}}\geq
\nu_{\overline{n-1}}.$ Hence $\sigma(\overline{n})=\overline{n}$ and $r_{\ast
}+2m_{\ast}=l$ that is $\sigma\in W_{n-1}.$ Moreover we have $\gamma
=\nu^{\prime}$ since $w\circ\nu^{\prime}=\gamma$ and $\nu,\gamma\in P_{n-1}.$
This proves that $s_{(\mu_{\overline{n}}+r+2m,\nu^{\prime})}$ can not be
obtained by applying the straightening law for Schur functions on
$s_{(\mu_{\overline{n}}+r_{\ast}+2m_{\ast},\gamma)}$ with $(\mu_{\overline{n}%
}+r_{\ast}+2m_{\ast},\gamma)\neq(\mu_{\overline{n}}+r+2m,\nu^{\prime}).$ Then
the corollary directly follows from (\ref{Morrisrec}) and Theorem
\ref{Th_hall_kostka}.
\end{proof}

\bigskip

Now suppose that $\nu=(p)_{n}=(p,0,...,0)\in P_{n}^{+}$.\ Then for any $\mu\in
P_{+}^{n}$ we must have%
\begin{equation}
K_{(p)_{n},\mu}(q)=\sum_{r+2m=l}q^{r+m}K_{(r)_{n-1},\mu^{\prime}}(q)
\label{for_re_row}%
\end{equation}
with $l=p-\mu_{\overline{n}}.$ This implies that $K_{(p)_{n},\mu}(q)$ may be
computed recursively. We are going to give an explicit formula for
$K_{(p)_{n},\mu}(q)$.\ The vertices of $B((p)_{n})_{\mu}=\{L\in B(p\Lambda
_{n-1}),$ $\mathrm{wt}(b)=\mu)\}$ are the words
\[
L=(n)^{k_{n}}\cdot\cdot\cdot(2)^{k_{2}}(1)^{k_{1}}(\overline{1})^{k_{\bar{1}}%
}(\overline{2})^{k_{\bar{2}}}\cdot\cdot\cdot(\overline{n})^{k_{\bar{n}}}%
\]
with $\mu_{\overline{i}}=k_{\overline{i}}-k_{i}$ for $i=1,...,n$ and
$k_{\overline{1}}+\cdot\cdot\cdot+k_{\overline{n}}+k_{1}+\cdot\cdot\cdot+k_{n}=p.$

\begin{proposition}
\label{prop_K_row}Let $p\geq1$ be an integer.\ For any $\mu\in P_{n}^{+}$ we
have%
\[
K_{(p)_{n},\mu}(q)=q^{f_{n}(\mu)}\sum_{L\in B((p)_{n})_{\mu}}q^{\theta_{n}(L)}%
\]
where $f_{n}(\mu)=\underset{i=1}{\overset{n}{\sum}}(n-i)\mu_{\overline{i}}$
and $\theta_{n}(L)=\underset{i=1}{\overset{n}{\sum}}(2\times
(n-i)+1)(k_{\overline{i}}-\mu_{\overline{i}}).$
\end{proposition}

\begin{proof}
We proceed by induction on $n.$ Suppose $n=1$ we have $f_{1}(\mu)=0.\;$We can
write $L=(1)^{p-k_{\bar{1}}}(\overline{1})^{k_{\bar{1}}}$ and $\mu
_{\overline{1}}=2k_{\overline{1}}-p.\;$Thus $\theta_{1}(L)=k_{\overline{1}%
}-\mu_{\overline{1}}=(p-\mu_{\overline{1}})/2$ and the proposition is true by
(\ref{Kostka_n=1}).

\noindent Now suppose the proposition true for $n-1.$ First note that
$f_{n}(\mu)=f_{n-1}(\mu^{\prime})+\underset{i=1}{\overset{n-1}{\sum}}%
\mu_{\overline{i}}.$ The set of vertices obtained by erasing the letters $n$
and $\overline{n}$ in $B((p)_{n})_{\mu}$ is the disjoint union of the
$B\left(  (r)_{n-1}\right)  _{\mu^{\prime}}$ with $r\in\{0,...,l=p-\mu
_{\overline{n}})$ since the number of letters $n$ or $\overline{n}$ belonging
to a vertex $L\in B((p)_{n})_{\mu}$ is a least equal to $\mu_{\overline{n}}.$
Its reflects the decomposition of $B((p)_{n})_{\mu}$ into its $U_{q}%
(sp_{2(n-1)})$-connected components. Consider $L\in B((p)_{n})_{\mu}$ and
denote by $L^{\prime}$ the vertex obtained by erasing all the letters $n$ and
$\overline{n}$ in $L.\;$Let $r$ be such that $L^{\prime}\in B\left(
(r)_{n-1}\right)  _{\mu^{\prime}}.\;$Then $r\in\{0,...,l\}$ and $l-r$ is even
since it is equal to the number of pairs $(n,\overline{n})$ erased in $L.\;$We
set $l-r=2m$.\ Then $k_{\overline{n}}=\mu_{\overline{n}}+m.$

\noindent We have
\[
\theta_{n}(L)=\theta_{n-1}(L^{\prime})+2\underset{i=1}{\overset{n-1}{\sum}%
}\left(  k_{\overline{i}}-\mu_{\overline{i}}\right)  +(k_{\overline{n}}%
-\mu_{\overline{n}}).
\]
From the equality $r=\underset{1\leq i\leq n-1}{\sum}\mu_{\overline{i}%
}+2\underset{1\leq i\leq n-1}{\sum}\left(  k_{\overline{i}}-\mu_{\overline{i}%
}\right)  $ we deduce%
\[
\theta_{n}(L)=\theta_{n-1}(L^{\prime})+r-\underset{i=1}{\overset{n-1}{\sum}%
}\mu_{\overline{i}}+m.
\]
Set%
\[
K=\sum_{L\in B((p)_{n})_{\mu}}q^{\theta_{n}(L)+f_{n}(\mu)}.
\]
Then by the above arguments%
\begin{multline*}
K=\sum_{r+2m=l}\sum_{L^{\prime}\in B\left(  (r)_{n-1}\right)  _{\mu^{\prime}}%
}q^{\theta_{n-1}(L^{\prime})+r-\underset{i=1}{\overset{n-1}{\sum}}%
\mu_{\overline{i}}+m+f_{n-1}(\mu^{\prime})+\underset{i=1}{\overset{n-1}{\sum}%
}\mu_{\overline{i}}}=\\
\sum_{r+2m=l}q^{r+m}\times\sum_{L^{\prime}\in B\left(  (r)_{n-1}\right)
_{\mu^{\prime}}}q^{\theta_{n-1}(L^{\prime})+f_{n-1}(\mu^{\prime})}.
\end{multline*}
Thus we obtain by the induction hypothesis%
\[
K=\sum_{r+2m=l}q^{r+m}K_{(r)_{n-1},\mu^{\prime}}(q).
\]
Finally $K=K_{(p)_{n},\mu}(q)$ by (\ref{for_re_row}).
\end{proof}

\begin{corollary}
\label{cor_h=2}Write $(1^{2})_{n}$ for the partition of length $n$ equal to
$(1,1,0,...,0)$. Then
\[
K_{(1^{2})_{n},0}(q)=\sum_{i=1}^{n-1}q^{2i}.
\]
\end{corollary}

\noindent To prove this corollary we need the more general lemma above

\begin{lemma}
\label{lem_dec_k_fund}Write $(1^{p})_{n}$ for the partition of length $n$
$(1,...1,0,...,0)$ with $p\geq2$ parts equal to $1.$ Then
\[
K_{(1^{p})_{n},0}(q)=(q-1)K_{\gamma^{p},0}(q)+qK_{(1^{p})_{n-1},0}%
(q)+qK_{(1^{p-2})_{n-1},0}(q)
\]
where $\gamma^{p}=(2,1,...,1,0,...,0)\in P_{n-1}^{+}$ contains $p-2$ parts
equal to $1.$
\end{lemma}

\begin{proof}
With $\mu=0$ formula (\ref{Morrisrec}) becomes%
\[
Q_{0}=\sum_{\gamma\in P_{n-1}^{+}}\sum_{r=0}^{+\infty}\sum_{m=0}^{+\infty
}q^{m+r}\sum_{\lambda\in B(\gamma)\otimes B((1)_{n-1}^{r})}K_{\lambda
,0}(q)s_{(r+2m,\gamma)}.
\]
By using the straightening law for Schur functions and Theorem
\ref{Th_hall_kostka} we have to find all the $(r+2m,\gamma)$ such that there
exits $\sigma\in W_{n-1}$ satisfying $\sigma\circ(r+2m,\gamma)=(1^{p})_{n}$
that is
\begin{equation}
\sigma(r+2m+n,\gamma_{\overline{n-1}}+n-1,...,\gamma_{\overline{1}%
}+1)=(n+1,...,n-p+2,n-p,...,1).\label{exp}%
\end{equation}
We have $r+2m+n\geq n$ hence $\sigma(\overline{n})\in\{\overline{n}%
,\overline{n-1}\}.$

\noindent$\mathrm{(i):}$ If $\sigma(\overline{n})=\overline{n}$ then $r=1$ and
$m=0.$ For $k\notin\{1,n\},$ $\gamma_{\overline{k}}+k>1$ thus $\sigma
(\overline{1})=\overline{1}.$ By a straightforward induction we obtain
$\sigma(\overline{k})=\overline{k}$ for $k\in\{1,...,n-p\}.\;$Moreover we have
$\gamma_{\overline{k}}+k\leq n$ for $k<n.\;$This implies that $\gamma
_{\overline{n-1}}\in\{0,1\}$ since $\gamma_{\overline{n-1}}+n-1\geq n-1$.\ We
can not have $\gamma_{\overline{n-1}}=0$ otherwise $\gamma_{\overline{k}}=0$
for any $k<n$ and the value $n$ in the left hand side of (\ref{exp}) is not
attained. Hence $\gamma_{\overline{n-1}}=1$ and $\sigma(\overline
{n-1})=\overline{n-1}.\;$By induction we can prove that $\gamma_{\overline
{n-1}}=\cdot\cdot\cdot=\gamma_{\overline{n-p+1}}=1$ and $\sigma(\overline
{k})=\overline{k}$ for $k\in\{n-1,...,n-p+1\}.$ It means that $\sigma=id,$
$R=1,$ $m=0$ and $\gamma=(1^{p-1})_{n-1}.$

\noindent$\mathrm{(ii):}$ If $\sigma(\overline{n})=\overline{n-1}$ then
$R=m=0.\;$By using similar arguments than above we obtain $\gamma
_{\overline{n-1}}=2,$ $\gamma_{\overline{n-2}}=\cdot\cdot\cdot=\gamma
_{\overline{n-p+1}}=1$ and $\gamma_{\overline{n-p}}=\cdot\cdot\cdot
\gamma_{\overline{1}}=0.$ It means that $\sigma=s_{n}$ and $\gamma=\gamma
^{p}.$ Note that $s_{(0,\gamma^{p})}=-s_{((1^{p})_{n}}$ since $s_{n}%
\circ(0,\gamma^{p})=(1^{p})_{n}$ and $l(s_{n})=1.$

\noindent Finally by Theorem \ref{Th_hall_kostka} we must have%
\begin{multline*}
K_{(1^{p})_{n},0}(q)=q\times\sum_{\lambda\in B((1^{p-1})_{n-1})\otimes
B((1)_{n-1})}K_{\lambda,0}(q)-K_{\gamma^{p},0}(q)=\\
(q-1)K_{\gamma^{p},0}(q)+qK_{(1^{p})_{n-1},0}+qK_{(1^{p-2})_{n-1},0}.
\end{multline*}
\end{proof}

\begin{proof}
(of Corollary \ref{cor_h=2}). We proceed by induction on $n.\;$For $n=2,$
$K_{(1^{2}),0}(q)=q^{2}$.\ Suppose the corollary true for $k<n.$ Then by
applying Lemma \ref{lem_dec_k_fund} we obtain%
\[
K_{(1^{2})_{n},0}(q)=(q-1)K_{(2)_{n-1},0}(q)+qK_{(1^{2})_{n-1},0}(q)+q.
\]
It follows from Proposition \ref{prop_K_row} that
\[
K_{(2)_{n-1},0}(q)=\sum_{i=1}^{n-1}q^{2i-1}.
\]
Thus%
\[
K_{(1^{2})_{n},0}(q)=(q-1)\sum_{i=1}^{n-1}q^{2i-1}+q\sum_{i=1}^{n-2}%
q^{2i}+q=\sum_{i=1}^{n-1}q^{2i}-\sum_{i=1}^{n-1}q^{2i-1}+\sum_{i=1}%
^{n-2}q^{2i+1}+q=\sum_{i=1}^{n-1}q^{2i}.
\]
\end{proof}

\noindent Note that we can not deduce an explicit formula for $K_{(1^{p}%
)_{n},0}(q)$ with $p>2$ from the recurrence formula of Lemma
\ref{lem_dec_k_fund} as we have do in Proposition \ref{prop_K_row} since we
have no explicit formula for $K_{\gamma^{p},0}(q)$ as soon as $p>2.$
Nevertheless we will give a conjectural general formula for $K_{(1^{p})_{n}%
,0}(q)$ in Section \ref{sec_cha}.

\section{Cyclage graphs for symplectic tableaux\label{sec_cy_graph}}

Given a symplectic tableau $T\in\mathbf{ST}(n)$, we can factorize
$\mathrm{w}(T)$ in a unique way by setting $\mathrm{w}(T)=xu$ where $u$ is a
word and $x$ is a letter.\ It is easy to verify that $u$ is also the reading
of a symplectic tableau, say $T_{\ast}\in\mathbf{ST}(n).$ The initial
cocyclage operation on $T$ consists in the insertion $x\rightarrow T_{\ast}.$
We are going to see that all the initial cocyclage operations are not relevant
for defining a charge.

\noindent It follows from Paragraph \ref{par_inser} that the tableau obtained
by cocycling a tableau $T\in\mathbf{ST}(n)$ does not belong to $\mathbf{ST}%
(n)$ in general but belongs to $\mathbf{ST}(n+1).\;$To overcome this problem
we are going to define our cocyclage operation directly on the complete
symplectic tableaux set $\mathbf{ST}=\underset{n\geq1}{\cup}\mathbf{ST}(n).$

\subsection{Cocyclage operation}

Set $\mathcal{C}_{\infty}=\underset{n\geq1}{\cup}\mathcal{C}_{n}$.\ Then
$\mathcal{C}_{\infty}$ is totally ordered by $\leq.\;$Given any $T\in
\mathbf{ST}$ there exists an integer $m\geq1$ such that $T\in\mathbf{ST}%
(m).\;$Recall that $d_{\overline{i}}$ is the number of letters $\overline{i}$
of $T$ minus the number of letters $i$.\ For any weight $\mu\in P_{n},$ we
will say that $T\in\mathbf{ST}$ is a tableau of weight $\mu$ if $T\in
\mathbf{ST}(m)$ with $m\geq n$, $d_{\overline{i}}=0$ for $i>n$ and
$d_{\overline{k}}=\mu_{\overline{k}}$ for $k=1,...,n.$ For any $\mu\in P_{n},$
the set of tableaux of weight $\mu$ is denoted $\mathbf{ST(}\mu\mathbf{).\;}%
$If $T\in\mathbf{ST(}\mu\mathbf{),}$ the number of letters $k$ with $k>n$
which belong to $T$ is equal to the number of letters $\overline{k}.$

\noindent Let $w\in\mathcal{C}_{\infty}^{\ast}$ and write $w=xu$ with $x$ a
letter and $u\in\mathcal{C}_{\infty}^{\ast}.$ The cocyclage shift $\xi$ is the
map defined on $\mathcal{C}_{\infty}^{\ast}$ by
\[
\xi(w)=ux.
\]

\begin{lemma}
\label{lem_kxi}For any $n\geq1$, $\sigma\in W_{n},$ and $w\in\mathcal{C}%
_{\infty}^{\ast}$, $\xi(\sigma(w))=\sigma(\xi(w)).$
\end{lemma}

\begin{proof}
The proof is analogous to that of Proposition 5.6.1 of \cite{loth}.
\end{proof}

\noindent Consider a symplectic tableau $T=C_{1}\cdot\cdot\cdot C_{r}%
\in\mathbf{ST}(m)$ with $r>1$.\ We will say that the cocyclage operation is
authorized for $T$ if there is no letter $y\in\mathcal{C}_{m}$ such that $y\in
C_{i}$ for any $i=1,...,r$ and $\overline{y}\notin T$. It means that the
cocyclage operation is not authorized for $T$ when there exits an integer
$p\in\{1,...,m\}$ such that $\left|  d_{\overline{p}}\right|  $ is equal to
$r.$ If the cocyclage operation is authorized for $T,$ we write $\mathrm{w}%
(T)=x\mathrm{w}(T_{\ast})$ where $T_{\ast}\in\mathbf{ST}(m)$ and
$x\in\mathcal{C}_{m}.\;$and we set%
\[
U(T)=x\rightarrow T_{\ast}.
\]

\noindent\textbf{Remarks:}

\noindent$\mathrm{(i)}\mathbf{:}$ $U(T)$ belongs to $\mathbf{ST}$ and
$\mathrm{wt}(U(T))=\mathrm{wt}(T).$ More precisely when $T\in\mathbf{ST}(m),$
$U(T)\in\mathbf{ST}(m)$ if the heights of the first columns of $T$ and $U(T)$
are equal, $U(T)\in\mathbf{ST}(m+1)$ otherwise.

\noindent$\mathrm{(ii)}\mathbf{:}$ If $\mathrm{wt}(T)=0$ then the cocyclage
operation is always authorized.

\noindent$\mathrm{(iii)}\mathbf{:}$ There is no cocyclage operation on the columns.

\begin{example}
Consider the tableaux $T_{1}=$%
\begin{tabular}
[c]{|l|ll}\hline
$\mathtt{\bar{4}}$ & $\mathtt{\bar{3}}$ & \multicolumn{1}{|l|}{$\mathtt{\bar
{2}}$}\\\hline
$\mathtt{\bar{2}}$ & $\mathtt{\bar{2}}$ & \multicolumn{1}{|l|}{$\mathtt{\bar
{1}}$}\\\hline
$\mathtt{2}$ &  & \\\cline{1-1}%
\end{tabular}
, $T_{2}=$%
\begin{tabular}
[c]{|l|ll}\hline
$\mathtt{\bar{4}}$ & $\mathtt{\bar{3}}$ & \multicolumn{1}{|l|}{$\mathtt{4}$%
}\\\hline
$\mathtt{\bar{2}}$ & $\mathtt{\bar{2}}$ & \multicolumn{1}{|l}{}\\\cline{1-2}%
$\mathtt{2}$ &  & \\\cline{1-1}%
\end{tabular}
and $T_{3}=$%
\begin{tabular}
[c]{|l|ll}\hline
$\mathtt{\bar{4}}$ & $\mathtt{\bar{3}}$ & \multicolumn{1}{|l|}{$\mathtt{\bar
{2}}$}\\\hline
$\mathtt{\bar{2}}$ & $\mathtt{\bar{2}}$ & \multicolumn{1}{|l|}{$\mathtt{\bar
{1}}$}\\\hline
$\mathtt{3}$ &  & \\\cline{1-1}%
\end{tabular}
. Then the cocyclage operation is authorized for $T_{1}$ and $T_{2}$ but not
in $T_{3}.\;$We obtain $U(T_{1})=$%
\begin{tabular}
[c]{|l|l|l}\hline
$\mathtt{\bar{4}}$ & $\mathtt{\bar{3}}$ & \multicolumn{1}{|l|}{$\mathtt{\bar
{1}}$}\\\hline
$\mathtt{\bar{3}}$ & $\mathtt{\bar{2}}$ & \\\cline{1-2}\cline{2-2}%
$\mathtt{\bar{2}}$ & $\mathtt{3}$ & \\\cline{1-2}%
\end{tabular}
, $U(T_{2})=%
\begin{tabular}
[c]{|l|l}\hline
$\mathtt{\bar{4}}$ & \multicolumn{1}{|l|}{$\mathtt{2}$}\\\hline
$\mathtt{\bar{2}}$ & \multicolumn{1}{|l|}{$\mathtt{3}$}\\\hline
$\mathtt{2}$ & \\\cline{1-1}%
$\mathtt{4}$ & \\\cline{1-1}%
\end{tabular}
$.
\end{example}

\begin{lemma}
\label{lem_U_com_sigma}Suppose $T\in\mathbf{ST}(m)$ and consider $\sigma\in
W_{m}.\;$Then the cocyclage operation is authorized for $T$ if and only if it
is authorized for $\sigma(T).\;$In this case we have $U(\sigma(T))=\sigma(U(T)).$
\end{lemma}

\begin{proof}
The lemma directly follows from Lemma \ref{lem_kxi} and (\ref{rq_p_com_w}).
\end{proof}

\subsection{Cyclage graphs}

We endow the set $\mathbf{ST}$ with a structure of graph by drawing an array
$T\rightarrow T^{\prime}$ if and only if the cocyclage operation is authorized
on $T$ and $U(T)=T^{\prime}.$ Write $\Gamma(T)$ for the connected component
containing $T.$ Let $t$ be the translation operation on letters of
$\mathcal{C}_{\infty}$ defined by $t(k)=k+1$ and $t(\overline{k}%
)=\overline{k+1}$ for $k\geq1$.\ We write $t(w)$ (resp.\ $t(T)$) for the word
(resp.\ the tableau) obtained by applying $t$ to each letter of $w\in
\mathcal{C}_{\infty}^{\ast}$ (resp.\ to each letter of $T\in\mathbf{ST)}$.

\begin{lemma}
\ \ \ \label{lem_shape}

\noindent$\mathrm{(i)}\mathbf{:}$ Suppose $T\in\mathbf{ST}(m)$.\ Then
$\Gamma(T)$ and $\Gamma(\sigma(T))$ are isomorphic for any $\sigma\in W_{m}$.

\noindent$\mathrm{(ii)}\mathbf{:}$ The cyclage graphs $\Gamma(T)$ and
$\Gamma(t(T))$ are isomorphic.

\noindent$\mathrm{(iii):}$ Suppose that $T_{1}\neq T_{2}\in\Gamma(T)$ are such
that $U(T_{1})=U(T_{2})=T.\;$Then $T_{1}$ and $T_{2}$ have different shapes.
\end{lemma}

\begin{proof}
\noindent Assertion $\mathrm{(i)}$ follows immediately from Lemma
\ref{lem_U_com_sigma}.

\noindent Let $w_{1}$ and $w_{2}$ be two words of $\mathcal{C}_{m}.$ Then
$w_{1}\equiv_{m}w_{2}$ if and only if $t(w_{1})\equiv_{m+1}t(w_{2})$.\ This
implies that $P(t(w))=t(P(w))$ for any word $w\in\mathcal{C}_{\infty}.$ Hence
$t$ commutes with $U$.\ Since $t$ is a bijection, it is also an isomorphism
between $\Gamma(T)$ and $\Gamma(t(T))$ which proves $\mathrm{(ii)}\mathbf{.}$

\noindent Suppose that $T_{1},T_{2}\in\mathbf{ST}(m).\;$Write $\mathrm{w}%
(T_{1})=x\mathrm{w}(R)$ and $\mathrm{w}(T_{2})=y\mathrm{w}(S)$ with $x,y$ two
letters and $R,S$ two symplectic tableaux.\ Then $\mathrm{w}(R)x\equiv
_{m+1}\mathrm{w}(S)y$ since $P(\mathrm{w}(R)x)=P(\mathrm{w}(S)y)=T.\;$Suppose
that $T_{1}$ and $T_{2}$ have the same shape.\ Then $R$ and $S$ have the same
shape $Y$.\ The highest weight vertices of the connected components of
$G_{m+1}$ containing $\mathrm{w}(R)x$ and $\mathrm{w}(S)y$ may be respectively
written $\mathrm{w}(Y_{0})x_{0}$ and $\mathrm{w}(Y_{0})y_{0}$ where $Y_{0}%
\in\mathbf{ST}(m+1)$ is the highest weight tableau of shape $Y.\;$The
congruence $\mathrm{w}(R)x\equiv_{m+1}\mathrm{w}(S)y$ implies the congruence
$\mathrm{w}(Y_{0})x_{0}\equiv_{m+1}\mathrm{w}(Y_{0})y_{0}.\;$Thus we must have
$\mathrm{wt}(\mathrm{w}(Y_{0})x_{0})=\mathrm{wt}(\mathrm{w}(Y_{0})x_{0}).\;$It
means that $x_{0}=y_{0}.\;$Hence $\mathrm{w}(R)x$ and $\mathrm{w}(S)y$ are
congruent and belong to the same connected component. This implies that
$\mathrm{w}(R)x=\mathrm{w}(S)y,$ thus $x\mathrm{w}(R)=y\mathrm{w}(S)$ and
$T_{1}=T_{2}.$ So $\mathrm{(iii)}$ is proved.
\end{proof}

\noindent Assertion $\mathrm{(i)}$ of the above lemma permits to restrict to
the cyclage graphs $\Gamma(T)$ with $T\in\mathbf{ST(}\mu\mathbf{)}$ and
$\mu\in P_{n}^{+}$. Suppose first that $\mu=0.$ We have seen that the
cocyclage operation is always authorized on symplectic tableaux of weight $0$
with at least two columns$.\;$So we can define from $T$ a sequence $(T_{n})$
of symplectic tableaux by setting $T_{0}=T$ and $T_{k+1}=U(T_{k})$ while $T$
is not a column.

\begin{proposition}
\label{lem_cyl_nul}The sequence $(T_{n})$ is finite without repetition and
there exists an integer $e$ such that $T_{e}$ is a column of weight $0.$
\end{proposition}

\noindent To prove this proposition we need two technical lemmas. Given two
words $w_{1},w_{2}\in\mathcal{C}_{\infty}^{\ast},$ write $w_{1}%
\vartriangleleft w_{2}$ if $w_{1}$ and $w_{2}$ can respectively be written
$w_{1}=u_{1}x_{1}v$ and $w_{2}=u_{2}x_{2}v$ where $u_{1},u_{2},v\in
\mathcal{C}_{\infty}^{\ast}$ and $x_{1},x_{2}\in\mathcal{C}_{\infty}$ verify
$x_{1}<x_{2}.$ It means that $\vartriangleleft$ is the inverse lexicographic
order on words of $\mathcal{C}_{\infty}^{\ast}.$ For any symplectic tableau
$T$ with $r>1$ columns, we denote by $N_{r}(T)$ the number of boxes belonging
to the $r-1$ rightmost columns of $T.$

\begin{lemma}
\label{lem_tech1}Consider $\mu\in P_{n}^{+}$ and $\tau\in\mathbf{ST(}%
\mu\mathbf{)}$ a tableau with $r>1\mathbf{\ }$columns.\ Let $T,T^{\prime}$ two
tableaux of $\Gamma\mathbf{(}\mu\mathbf{)}$ such that $T=U^{(i)}(\tau)$ and
$T^{\prime}=U^{(i+1)}(\tau)$ with $i\geq0$ an integer. Then the following
assertions hold.

\begin{enumerate}
\item $T$ and $T^{\prime}$ contains at most $r+1$ columns.

\item  Suppose that $T$ contains $r$ or $r+1$ columns and $N_{r}%
(T)=N_{r}(T^{\prime})$.\ Then only one of the following situations can happen:

\noindent$\mathrm{(i):}$ $T$ and $T^{\prime}$ contain $r$ columns and their
$r$-th columns have the same height.

\noindent$\mathrm{(ii):}$ $T$ contains $r$ columns and $T^{\prime}$ contains
$r+1$ columns.

\noindent$\mathrm{(iii):}$ $T$ and $T^{\prime}$ contains $r+1$ columns.

\noindent$\mathrm{(iv):}$ $T$ contains $r+1$ columns, $T$ contains $r$ columns
and the height of the last column $C_{r}^{\prime}$ of $T^{\prime}$ is equal to
$h(C_{r})+1$.

Moreover in each case we can write $\mathrm{w}(T)=x_{\ast}\mathrm{w}(T_{\ast
})$ and $\mathrm{w}(T^{\prime})=x_{\ast}^{\prime}\mathrm{w}(T_{\ast}^{\prime
})$ with $\mathrm{w}(T_{\ast}^{\prime})\vartriangleleft\mathrm{w}(T_{\ast}).$

\item  Suppose that $T$ contains $r$ or $r+1$ columns and $N_{r}(T)\neq
N_{r}(T^{\prime})$.\ Then $N_{r}(T)>N_{r}(T^{\prime})$.
\end{enumerate}
\end{lemma}

\begin{proof}
$1:$ Let $j$ be an integer such that $U^{(j)}(\tau)$ is defined and contains
$r+1$ columns.\ Then the height of its last column is equal to $1.$ The
assertion follows immediately.

$2:$ In case $\mathrm{(i)}$ the box which is added to $T_{\ast}$ during the
insertion $x_{\ast}\rightarrow T_{\ast}$ appears on the bottom of the last
column $C_{r,\ast}$ (eventually empty) of $T_{\ast}.\;$This insertion can be
written $x_{\ast}\rightarrow T_{\ast}=(x_{\ast}\rightarrow C_{1}\cdot
\cdot\cdot C_{r-1})C_{r,\ast}=(C_{1}^{\prime}\cdot\cdot\cdot C_{r-1}^{\prime
})(y\rightarrow C_{r,\ast})$ that is, $x_{\ast}$ is first inserted in the
sub-tableau composed of the $r-1$ leftmost columns of $T_{\ast}$ which gives a
new tableau $C_{1}^{\prime}\cdot\cdot\cdot C_{r-1}^{\prime}$ and a letter
$y$.\ This letter is then inserted on the bottom of $C_{r,\ast}.$ Suppose that
there exists an integer $i\in\{1,...,r-1\}$ such that $C_{i}^{\prime}\neq
C_{i}.\;$Then if we choose $i$ minimal we have $\mathrm{w}(C_{i}^{\prime
})\vartriangleleft\mathrm{w}(C_{i})$ by (\ref{C'infC}) and finally
$\mathrm{w}(T_{\ast}^{\prime})\vartriangleleft\mathrm{w}(T_{\ast})$. Now if
$C_{i}^{\prime}=C_{i}$ for $i=1,...,r-1$ we have $y=x_{\ast}$ and $x_{\ast}\in
C_{i},$ $\overline{x}_{\ast}\notin C_{i}$ for any $i=1,...,r-1$. So the letter
$x_{\ast}$ belongs to all the columns of $T$.\ Then $x_{\ast}$ is a barred
letter since $\mu\in P_{n}^{+}$ and $r>1.\;$Moreover $\overline{x}_{\ast
}\notin C_{r}^{\ast}$ for $x_{\ast}\rightarrow C_{r}^{\ast}$ is a
column.\ Thus $\overline{x}_{\ast}\notin T.$ This contradicts the fact that
the cocyclage operation is authorized for $T.$

\noindent In case $\mathrm{(ii)}$ a new column of height $1$ is added to the
shape of $T_{\ast}.$ The insertion can be written $x_{\ast}\rightarrow
T_{\ast}=(x_{\ast}\rightarrow C_{1}\cdot\cdot\cdot C_{r-1})C_{r,\ast}%
=(C_{1}^{\prime}\cdot\cdot\cdot C_{r-1}^{\prime})(y\rightarrow C_{r,\ast
})=C_{1}^{\prime}\cdot\cdot\cdot C_{r-1}^{\prime}C_{r,\ast}^{\prime}%
\begin{tabular}
[c]{|l|}\hline
$x_{\ast}^{\prime}$\\\hline
\end{tabular}
$ that is, $x_{\ast}$ is inserted in the sub-tableau composed of the $r-1$
leftmost columns of $T_{\ast}$ which gives a new tableau $C_{1}^{\prime}%
\cdot\cdot\cdot C_{r-1}^{\prime}$ and a letter $y$.\ This letter is then
inserted in $C_{r,\ast}$ which gives the column $C_{r,\ast}^{\prime}$ and the
letter $x_{\ast}^{\prime}.$ If $y\neq x_{\ast},$ we terminate as in case
$\mathrm{(i).}$ Otherwise we have $C_{i}^{\prime}=C_{i}$ for any
$i=1,...,r-1.$ We can not have $x_{\ast}^{\prime}=x_{\ast}$ since it would
imply that $x_{\ast}\in C_{r,\ast}$ which is impossible since $C_{r}$ can not
contain two letters $x_{\ast}.$ Thus $x_{\ast}^{\prime}>x_{\ast}$,
$\mathrm{w}(C_{r,\ast}^{\prime})\vartriangleleft\mathrm{w}(C_{r,\ast})$ and
finally $\mathrm{w}(T_{\ast}^{\prime})\vartriangleleft\mathrm{w}(T_{\ast}).$

\noindent Case $\mathrm{(iii)}$ is similar to case $\mathrm{(i)}$ with
$h(C_{r})=1.$

\noindent In case $\mathrm{(iv)}$ $C_{r+1}$ contains only the letter $x_{\ast
}$ and a new box appears on the bottom of the column $C_{r}$ of $T_{\ast}$
during the insertion $x_{\ast}\rightarrow T_{\ast}$. The insertion can be
written $x_{\ast}\rightarrow T_{\ast}=x_{\ast}\rightarrow(C_{1}\cdot\cdot\cdot
C_{r})=(C_{1}^{\prime}\cdot\cdot\cdot C_{r-1}^{\prime})(y\rightarrow
C_{r,\ast})=C_{1}^{\prime}\cdot\cdot\cdot C_{r-1}^{\prime}C_{r,\ast}^{\prime
},$ that is $x_{\ast}$ is inserted in the sub-tableau composed of the $r-1$
leftmost columns of $T_{\ast}$ which gives the tableau $(C_{1}^{\prime}%
\cdot\cdot\cdot C_{r-1}^{\prime})$ and the letter $y.\;$This letter is then
inserted on the bottom of $C_{r,\ast}$ which gives the column $C_{r,\ast
}^{\prime}.$ Suppose that $y=x_{\ast}.\;$We must have $x_{\ast}\in C_{i},$ and
$\overline{x}_{\ast}\notin C_{i}$ for any $i=1,...,r-1$. Then $x_{\ast
}=\overline{q}$ with $q\geq1$ that is, is a barred letter as in $\mathrm{(i)}%
.\;$Moreover $\overline{x}_{\ast},x_{\ast}\notin C_{r}$ because $x_{\ast
}\rightarrow C_{r}$ is a column.\ Thus $d_{\overline{q}}(\mathrm{w}(T))=r.$
Now $T$ contains $r+1$ columns, hence $T\neq\tau$.\ Let $j$ minimal such that
$R=U^{(i-j)}(\tau)$ contains $r$ columns. Then $d_{\overline{q}}%
(\mathrm{w}(R))=d_{\overline{q}}(\mathrm{w}(T))=r.$ Thus the cocyclage
operation in not authorized in $R$ and we obtain a contradiction. It means
that $y\neq x_{\ast}.\;$So we can terminate as in case $\mathrm{(i).}$

$3:$ It is clear from the definition of the cocyclage operation.
\end{proof}

\bigskip

\begin{lemma}
\label{lem_tec_2}Let $T=C_{1}\cdot\cdot\cdot C_{p}\in\mathbf{ST(}0\mathbf{)}$
with $r>1$ columns.\ Then there exists an integer $k$ such that $T_{k}$ has at
most $r-1$ columns. Moreover if $k$ is minimal the sequence $T_{0},...,T_{k}$
is without repetition.
\end{lemma}

\begin{proof}
Let $m\geq1$ be an integer such that $T_{0},...,T_{m}$ have $r$ or $r+1$
columns and $N_{r}(T_{i})=N_{r}(T)$ for any $i=1,...,m.$ Then by Lemma
\ref{lem_tech1} we can write $\mathrm{w}(T_{i})=x_{i,\ast}\mathrm{w}%
(T_{i,\ast})$, $i=0,...,m$ with
\begin{equation}
\mathrm{w}(T_{m,\ast})\vartriangleleft\cdot\cdot\cdot\vartriangleleft
\mathrm{w}(T_{1,\ast})\vartriangleleft\mathrm{w}(T_{0,\ast})\text{.}
\label{seq_mono}%
\end{equation}
This implies that the sequence $T_{0},...,T_{m}$ is without repetition. Now
suppose that $T\in\mathbf{ST(}n\mathbf{).}$ Then $T_{i}\in\mathbf{ST(}%
n\mathbf{)}$ for any $i=0,...,m$.\ Indeed the height of the first column of
any tableau $T_{i},$ $i=0,...,m$ is always equal to that of $T_{0}$ since
$N_{r}(T_{i})=N_{r}(T)$ (see Remark $\mathrm{(i)}$ after Lemma \ref{lem_kxi}).
Denote by $p$ the number of boxes in $T_{0}.\;$Since the number of symplectic
tableaux with $p$ boxes belonging to $\mathbf{ST}(n)$ is finite there exits an
integer $s_{1}$ minimal such that $N_{r}(T_{s_{1}})=N_{r}(T)+1$ or $T_{s_{1}}$
has at most $r-1$ columns.\ Then the sequence $T_{0},...,T_{s_{1}}$ is without
repetition. If $T_{s_{1}}$ has at most $r-1$ columns we take $k=s_{1}%
.\;$Otherwise $T_{s_{1}}$ has $r$ or $r+1$ columns and we can obtain similarly
starting from $T_{s_{1}}$ an integer $s_{2}$ minimal such that $N_{p}%
(T_{s_{2}})=N_{p}(T_{s_{1}})+1$ or $T_{s_{2}}$ has at most $r-1$ columns.\ The
sequence $T_{s_{1}},...,T_{s_{2}}$ is without repetition.\ Then the sequence
$T_{0},...,T_{s_{2}}$ is also without repetition. Indeed a tableau $T_{i}$
with $i\in\{0,...,s_{1}-1\}$ can not be equal to a tableau $T_{j}$ with
$j\in\{s_{1},...,s_{2}-1\}$ since $N_{r}(T_{i})\neq N_{r}(T_{j}).\;$By
induction we can construct $T_{s_{j+1}}$ from $T_{s_{j}}$ while $T_{s_{j}}$
has $r$ or $r+1$ columns, such that $N_{r}(T_{s_{j+1}})=N_{r}(T_{s_{j}})+1$
and the sequence $T_{0},...,T_{s_{j+1}}$ is without repetition. The procedure
terminates since the number of boxes belonging to the columns $r$ and $r+1$
decreases by $1$ to each step. So the lemma is proved.
\end{proof}

\bigskip

\begin{proof}
(of Proposition \ref{lem_cyl_nul})

\noindent Let $r>1$ be the number of columns of $T.\;$By Lemma \ref{lem_tec_2}%
, we can obtain from $T=T_{0}$ a tableau $T_{k_{1}}$ with at most $r-1$
columns and such that the sequence $T_{0},...,T_{k_{1}}$ is without
repetition.\ If $r-1>1$ we can obtain a tableau $T_{k_{2}}$ from $T_{k_{1}}$
with at most $r-2$ columns and such that the sequence $T_{0},...,T_{k_{2}}$ is
without repetition. We can define $T_{k_{s+1}}$ from $T_{k_{s}}$ while $r-s>1$
such that the sequence $T_{0},...,T_{k_{s+1}}$ is without repetition. It is
clear that the procedure terminates when $T_{k_{s}}=T_{e}$ is a column of
weight $0.$
\end{proof}

\bigskip

It follows from Proposition \ref{lem_cyl_nul} that $\mathrm{wt}(T_{1}%
)=\mathrm{wt}(T_{2})\not \Longrightarrow\Gamma(T_{1})=\Gamma(T_{2})$ in
general. For example all the columns of weight $0$ occur in different
connected components.

\noindent We give below $\Gamma\left(
\begin{tabular}
[c]{|l|l|l|l|}\hline
$\mathtt{\bar{1}}$ & $\mathtt{\bar{1}}$ & $\mathtt{1}$ & $\mathtt{1}$\\\hline
\end{tabular}
\right)  ,$ $\Gamma\left(
\begin{tabular}
[c]{|l|l|l|l|l|}\hline
$\mathtt{\bar{3}}$ & $\mathtt{\bar{3}}$ & $\mathtt{\bar{2}}$ & $\mathtt{\bar
{1}}$ & $\mathtt{1}$\\\hline
\end{tabular}
\right)  ,$ $\Gamma\left(
\begin{tabular}
[c]{|l|l|}\hline
$\mathtt{\bar{3}}$ & $\mathtt{1}$\\\hline
$\mathtt{\bar{2}}$ & $\mathtt{2}$\\\hline
$\mathtt{\bar{1}}$ & $\mathtt{3}$\\\hline
\end{tabular}
\right)  $ and $\Gamma\left(
\begin{tabular}
[c]{|l|l|}\hline
$\mathtt{\bar{2}}$ & $\mathtt{1}$\\\hline
$\mathtt{\bar{1}}$ & $\mathtt{2}$\\\hline
\end{tabular}
\right)  $%

\[%
\begin{tabular}
[c]{ccc}%
\begin{tabular}
[c]{|l|}\hline
$\mathtt{\bar{3}}$\\\hline
$\mathtt{\bar{1}}$\\\hline
$\mathtt{1}$\\\hline
$\mathtt{3}$\\\hline
\end{tabular}
$\vspace{0.15cm}$ &  & \\
$\uparrow\vspace{0.15cm}$ &  & \\%
\begin{tabular}
[c]{|l|l}\hline
$\mathtt{\bar{3}}$ & \multicolumn{1}{|l|}{$\mathtt{3}$}\\\hline
$\mathtt{\bar{1}}$ & \\\cline{1-1}%
$\mathtt{1}$ & \\\cline{1-1}%
\end{tabular}
$\vspace{0.15cm}$ &  & \\
$\uparrow\vspace{0.15cm}$ & $\nwarrow$ & \\%
\begin{tabular}
[c]{|l|l}\hline
$\mathtt{\bar{2}}$ & \multicolumn{1}{|l|}{$\mathtt{1}$}\\\hline
$\mathtt{\bar{1}}$ & \\\cline{1-1}%
$\mathtt{2}$ & \\\cline{1-1}%
\end{tabular}
$\vspace{0.15cm}$ &  & $%
\begin{tabular}
[c]{|l|l|}\hline
$\mathtt{\bar{3}}$ & $\mathtt{1}$\\\hline
$\mathtt{\bar{1}}$ & $\mathtt{3}$\\\hline
\end{tabular}
$\\
$\uparrow\vspace{0.15cm}$ & $\nwarrow$ & \\%
\begin{tabular}
[c]{|l|ll}\hline
$\mathtt{\bar{2}}$ & $\mathtt{1}$ & \multicolumn{1}{|l|}{$\mathtt{2}$}\\\hline
$\mathtt{\bar{1}}$ &  & \\\cline{1-1}%
\end{tabular}
$\vspace{0.15cm}$ &  &
\begin{tabular}
[c]{|l|l}\hline
$\mathtt{\bar{2}}$ & \multicolumn{1}{|l|}{$\mathtt{\bar{1}}$}\\\hline
$\mathtt{1}$ & \\\cline{1-1}%
$\mathtt{2}$ & \\\cline{1-1}%
\end{tabular}
\\
$\uparrow\vspace{0.15cm}$ &  & $\uparrow$\\
$%
\begin{tabular}
[c]{|l|l|}\hline
$\mathtt{\bar{2}}$ & $\mathtt{\bar{1}}$\\\hline
$\mathtt{1}$ & $\mathtt{2}$\\\hline
\end{tabular}
\vspace{0.15cm}$ &  &
\begin{tabular}
[c]{|l|ll}\hline
$\mathtt{\bar{2}}$ & $\mathtt{\bar{1}}$ & \multicolumn{1}{|l|}{$\mathtt{2}$%
}\\\hline
$\mathtt{1}$ &  & \\\cline{1-1}%
\end{tabular}
\\
$\uparrow\vspace{0.15cm}$ &  & \\%
\begin{tabular}
[c]{|l|ll}\hline
$\mathtt{\bar{1}}$ & $\mathtt{\bar{1}}$ & \multicolumn{1}{|l|}{$\mathtt{1}$%
}\\\hline
$\mathtt{1}$ &  & \\\cline{1-1}%
\end{tabular}
$\vspace{0.15cm}$ &  & \\
$\uparrow\vspace{0.15cm}$ &  & \\
$%
\begin{tabular}
[c]{|l|l|l|l|}\hline
$\mathtt{\bar{1}}$ & $\mathtt{\bar{1}}$ & $\mathtt{1}$ & $\mathtt{1}$\\\hline
\end{tabular}
$ &  &
\end{tabular}
\]%

\begin{equation}%
\begin{tabular}
[c]{ccccc}%
$\vspace{0.15cm}$ &  &
\begin{tabular}
[c]{|l|l}\hline
$\mathtt{\bar{3}}$ & \multicolumn{1}{|l|}{$\mathtt{\bar{3}}$}\\\hline
$\mathtt{\bar{2}}$ & \multicolumn{1}{|l|}{$\mathtt{\bar{1}}$}\\\hline
$\mathtt{1}$ & \\\cline{1-1}%
\end{tabular}
&  & \\
$\vspace{0.15cm}$ & $\nearrow$ &  & $\nwarrow$ & \\
$\vspace{0.15cm}$%
\begin{tabular}
[c]{|l|ll}\hline
$\mathtt{\bar{3}}$ & $\mathtt{\bar{3}}$ & \multicolumn{1}{|l|}{$\mathtt{\bar
{2}}$}\\\hline
$\mathtt{\bar{1}}$ &  & \\\cline{1-1}%
$\mathtt{1}$ &  & \\\cline{1-1}%
\end{tabular}
&  &  &  & $%
\begin{tabular}
[c]{|l|l|l}\hline
$\mathtt{\bar{3}}$ & $\mathtt{\bar{3}}$ & \multicolumn{1}{|l|}{$\mathtt{1}$%
}\\\hline
$\mathtt{\bar{2}}$ & $\mathtt{\bar{1}}$ & \\\cline{1-2}%
\end{tabular}
$\\
$\vspace{0.15cm}\uparrow$ &  &  &  & $\uparrow$\\
$\vspace{0.15cm}%
\begin{tabular}
[c]{|l|lll}\hline
$\mathtt{\bar{3}}$ & $\mathtt{\bar{3}}$ & \multicolumn{1}{|l}{$\mathtt{\bar
{2}}$} & \multicolumn{1}{|l|}{$\mathtt{1}$}\\\hline
$\mathtt{\bar{1}}$ &  &  & \\\cline{1-1}%
\end{tabular}
$ &  &  &  & $%
\begin{tabular}
[c]{|l|l|l}\hline
$\mathtt{\bar{3}}$ & $\mathtt{\bar{3}}$ & \multicolumn{1}{|l|}{$\mathtt{\bar
{2}}$}\\\hline
$\mathtt{\bar{1}}$ & $\mathtt{1}$ & \\\cline{1-2}%
\end{tabular}
$\\
&  &  &  & $\vspace{0.15cm}\uparrow$\\
&  &  &  & $\vspace{0.15cm}%
\begin{tabular}
[c]{|l|lll}\hline
$\mathtt{\bar{3}}$ & $\mathtt{\bar{3}}$ & \multicolumn{1}{|l}{$\mathtt{\bar
{2}}$} & \multicolumn{1}{|l|}{$\mathtt{\bar{1}}$}\\\hline
$\mathtt{1}$ &  &  & \\\cline{1-1}%
\end{tabular}
$\\
&  &  &  & $\vspace{0.15cm}\uparrow$\\
&  &  &  &
\begin{tabular}
[c]{|l|l|l|l|l|}\hline
$\mathtt{\bar{3}}$ & $\mathtt{\bar{3}}$ & $\mathtt{\bar{2}}$ & $\mathtt{\bar
{1}}$ & $\mathtt{1}$\\\hline
\end{tabular}
\end{tabular}
\qquad%
\begin{tabular}
[c]{c}%
$\vspace{0.15cm}%
\begin{tabular}
[c]{|l|}\hline
$\mathtt{\bar{3}}$\\\hline
$\mathtt{\bar{2}}$\\\hline
$\mathtt{\bar{1}}$\\\hline
$\mathtt{1}$\\\hline
$\mathtt{2}$\\\hline
$\mathtt{3}$\\\hline
\end{tabular}
$\\
$\vspace{0.15cm}\uparrow$\\
$\vspace{0.15cm}%
\begin{tabular}
[c]{|l|l}\hline
$\mathtt{\bar{3}}$ & \multicolumn{1}{|l|}{$\mathtt{3}$}\\\hline
$\mathtt{\bar{2}}$ & \\\cline{1-1}%
$\mathtt{\bar{1}}$ & \\\cline{1-1}%
$\mathtt{1}$ & \\\cline{1-1}%
$\mathtt{2}$ & \\\cline{1-1}%
\end{tabular}
$\\
$\vspace{0.15cm}\uparrow$\\
$\vspace{0.15cm}%
\begin{tabular}
[c]{|l|l}\hline
$\mathtt{\bar{3}}$ & \multicolumn{1}{|l|}{$\mathtt{2}$}\\\hline
$\mathtt{\bar{2}}$ & \multicolumn{1}{|l|}{$\mathtt{3}$}\\\hline
$\mathtt{\bar{1}}$ & \\\cline{1-1}%
$\mathtt{1}$ & \\\cline{1-1}%
\end{tabular}
$\\
$\vspace{0.15cm}\uparrow$\\
$%
\begin{tabular}
[c]{|l|l|}\hline
$\mathtt{\bar{3}}$ & $\mathtt{1}$\\\hline
$\mathtt{\bar{2}}$ & $\mathtt{2}$\\\hline
$\mathtt{\bar{1}}$ & $\mathtt{3}$\\\hline
\end{tabular}
$%
\end{tabular}
\ \qquad\ \
\begin{tabular}
[c]{c}%
$\vspace{0.15cm}$%
\begin{tabular}
[c]{|l|}\hline
$\mathtt{\bar{2}}$\\\hline
$\mathtt{\bar{1}}$\\\hline
$\mathtt{1}$\\\hline
$\mathtt{2}$\\\hline
\end{tabular}
\\
$\vspace{0.15cm}\uparrow$\\
$\vspace{0.15cm}$%
\begin{tabular}
[c]{|l|l}\hline
$\mathtt{\bar{2}}$ & \multicolumn{1}{|l|}{$\mathtt{2}$}\\\hline
$\mathtt{\bar{1}}$ & \\\cline{1-1}%
$\mathtt{1}$ & \\\cline{1-1}%
\end{tabular}
\\
$\vspace{0.15cm}\uparrow$\\
$%
\begin{tabular}
[c]{|l|l|}\hline
$\mathtt{\bar{2}}$ & $\mathtt{1}$\\\hline
$\mathtt{\bar{1}}$ & $\mathtt{2}$\\\hline
\end{tabular}
$%
\end{tabular}
\label{graph(2,1,0)}%
\end{equation}

\noindent\textbf{Remarks:}

\noindent$\mathrm{(i)}\mathbf{:}$ Given $T^{\prime}\in\mathbf{ST,}$ it is
possible to find the tableaux $T$ (if there is any) such that $U(T)=T^{\prime
}.$ To do this we find all the pairs $(x,T_{\ast})$ obtained by applying the
reverse insertion algorithm on the outside corners of $T^{\prime}.$ By
definition of $U,$ the tableaux $T$ are precisely those which verify
$\mathrm{w}(T)=x\mathrm{w}(T_{\ast})$ for a pair $(x,T_{\ast})$. They are
determined by the pairs $(x,T_{\ast})$ for which $x\mathrm{w}(T_{\ast})$ is
the reading of a symplectic tableau. For example $T^{\prime}=%
\begin{tabular}
[c]{|l|l|}\hline
$\mathtt{\bar{2}}$ & $\mathtt{1}$\\\hline
$\mathtt{\bar{1}}$ & $\mathtt{2}$\\\hline
\end{tabular}
$ has only one outside corner which gives $x=\overline{1}$ and $T_{\ast}=%
\begin{tabular}
[c]{|l|l}\hline
$\mathtt{\bar{1}}$ & \multicolumn{1}{|l|}{$\mathtt{1}$}\\\hline
$\mathtt{1}$ & \\\cline{1-1}%
\end{tabular}
.$ There is no tableau $T$ such that $U(T)=T^{\prime}$ since $\overline
{1}(1\overline{1}1)$ is not the reading of a symplectic tableau.

\noindent$\mathrm{(ii)}\mathbf{:}$ In the definition of $U(T)$ we have
restricted to the initial cocyclages. For type $A$ the cyclage graphs take
also into account non initial cyclages.\ When $T$ is of dominant evaluation,
they are obtained by considering all the factorizations $\mathrm{w}(T)\equiv
y\mathrm{w}(Y)$ in the plactic monoid with $y\neq1$ a letter and $Y$ a
semi-standard tableau$.$ The use of non initial cocyclages with symplectic
tableaux is problematic because it can make appear loops in the cyclages
graphs. For example consider $Z=$%
\begin{tabular}
[c]{|l|l}\hline
$\mathtt{\bar{3}}$ & \multicolumn{1}{|l|}{$\mathtt{3}$}\\\hline
$\mathtt{\bar{1}}$ & \\\cline{1-1}%
$\mathtt{1}$ & \\\cline{1-1}%
\end{tabular}
.\ Then $\mathrm{w}(Z)=(3\overline{3}\ \overline{1})1\equiv_{3}\overline
{2}(2\overline{1}1)\equiv_{3}\overline{2}\ \overline{1}21.$ Now if we compute
$P(\xi(\overline{2}\ \overline{1}21))$ we obtain $Z^{\prime}=P(\overline
{1}21\overline{2})=$%
\begin{tabular}
[c]{|l|ll}\hline
$\mathtt{\bar{2}}$ & $\mathtt{\bar{1}}$ & \multicolumn{1}{|l|}{$\mathtt{2}$%
}\\\hline
$\mathtt{1}$ &  & \\\cline{1-1}%
\end{tabular}
. So we have a loop since $U^{(3)}(Z^{\prime})=Z$ (see the cyclage graph
$\Gamma\left(
\begin{tabular}
[c]{|l|l|l|l|}\hline
$\mathtt{\bar{1}}$ & $\mathtt{\bar{1}}$ & $\mathtt{1}$ & $\mathtt{1}$\\\hline
\end{tabular}
\right)  $ above).

\noindent$\mathrm{(iii)}\mathbf{:}$ For type $A,$ every semi-standard tableau
belongs to the cyclage graph containing a row tableau. By considering
$\Gamma\left(
\begin{tabular}
[c]{|l|l|}\hline
$\mathtt{\bar{2}}$ & $\mathtt{1}$\\\hline
$\mathtt{\bar{1}}$ & $\mathtt{2}$\\\hline
\end{tabular}
\right)  $ we see that such a property is false with the symplectic tableaux
even if we consider non initial cyclages. This explains why we have to
consider the cocyclage operation and not the cyclage one.

\subsection{Reduction operations}

\noindent Consider $T\in\mathbf{ST}(\mu)$ with $\mu\in P_{n}^{+}.\;$If the
cocyclage operation is not authorized for $T,$ then $T$ do not contain any
letters $n$.\ Indeed there exists $p\in\{1,...n\}$ such that $\mu
_{\overline{p}}$ is equal to the number of columns of $T.\;$Thus
$\mu_{\overline{n}}=\mu_{\overline{p}}$ since $\mu_{\overline{n}}\geq
\mu_{\overline{p}}$.\ So each column of $T$ contains a letter $\overline{n}$
and no letter $n.$ Let $T_{\$}$ be the tableau obtained first by erasing the
letters $\overline{n}$ in $T$ next by applying $t$ to the letters $x\in T$
such that $\overline{n}<x<n$.\ It is easy to verify that $T_{\$}\in
\mathbf{ST}(\mu^{\prime})$ with $\mu^{\prime}=(\mu_{\overline{n-1}}%
,...,\mu_{\overline{1}},0)\in P_{n}^{+}$.\ Now if the cocyclage operation is
not authorized on $T_{\$}$, we can compute $(T_{\$})_{\$}$ and so on until
obtain a symplectic tableau $\widehat{T}$ which is either a column (eventually
empty) either a symplectic tableau for which the cocyclage operation is
authorized. We will say that $\widehat{T}$ is obtained by reduction operations
from $T.$ By convention we set $\widehat{T}=T$\ if the cocyclage operation is
already defined for the symplectic tableau $T.$

\noindent\textbf{Remark: }When a reduction operation is done in $T\in
\mathbf{ST}(\mu)$ with $\mu\in P_{n}^{+},$ $\mu_{\overline{n}}$ is equal to
the first part of the shape of $T.\;$To define a charge statistic on
symplectic tableaux related to Kostka-Foulkes polynomials, it seems natural by
Lemma \ref{lem_lamb1=mu1} to impose that the charges associated to $T$ and
$\widehat{T}$ should be equal as we will do in Section \ref{sec_cha}.

\bigskip

From $T\in\mathbf{ST}(\mu)$ we can compute a sequence of symplectic tableaux
by setting $T_{0}=T$ and
\[
T_{k+1}=U(\widehat{T}_{k})
\]
while $\widehat{T}_{k}$ is not a column.

\begin{proposition}
\label{propcy2}The sequence $(T_{n})$ is finite without repetition and the
last symplectic tableau obtained is a column of weight $0$ (eventually empty).
\end{proposition}

\begin{proof}
Suppose first that there is a loop in the sequence $(T_{n})$ that is, there
exits two integers $k$ and $s$ such that $T_{k}=T_{k+s}.$ Then $T_{i}%
=\widehat{T}_{i}$ for any $i=k,...,k+s-1.$ Choose $p\in\{k,..,k+s-1\}$ such
that the number of columns of $T_{p}$ is minimal among all the tableaux
$T_{i},$ $i=k,...,k+s-1.\;$Denote by $r$ the number of columns of $T_{p}.$
Then by assertion $1$ of Lemma \ref{lem_tech1}$,T_{k},...,T_{k+s-1}$ contain
$r$ or $r+1$ columns since every $T_{i},$ $i=k,...,k+s-1$ can be obtained by
cocyclage operations from $T_{p}$. We must have $N_{r}(T_{k})\leq\cdot
\cdot\cdot\leq N_{r}(T_{k+s}).\;$This implies that $N_{r}(T_{k})=\cdot
\cdot\cdot=N_{r}(T_{k+s-1})$ for $T_{k+s}=T_{k}.$ Then by assertion $2$ of
Lemma \ref{lem_tech1} we can write $\mathrm{w}(T_{i})=x_{i,\ast}%
\mathrm{w}(T_{i,\ast})$,$i=k,...,k+s$ with
\[
\mathrm{w}(T_{k+s,\ast})\vartriangleleft\cdot\cdot\cdot\vartriangleleft
\mathrm{w}(T_{s,\ast}).
\]
We obtain a contradiction since $\mathrm{w}(T_{k+s,\ast})=\mathrm{w}%
(T_{s,\ast}).$ It means that there is no loop in the sequence $(T_{n})$. Hence
this sequence is without repetition.

\noindent Now suppose that this sequence is infinite.\ Then there exits an
integer $a$ such that the sequence $(T_{n+a})_{n\geq0}$ is infinite without
reduction operation. In the proof of Lemma \ref{lem_tec_2} the hypothesis
$\mu=0$ is only used to assure that the sequence of the cocycled tableaux is
defined. It means that this lemma is yet true for the sequence $(T_{n+a}%
)_{n\geq0}.$ Thus we can define by induction as in proof of Proposition
\ref{lem_cyl_nul} an infinite sequence of tableaux $(T_{v_{j}})_{j\geq0}$ such
that $T_{j_{0}}=T_{n+a}$ and for any $j,$ $T_{v_{j+1}}$ has one column less
than $T_{v_{j}}.$ We derive a contradiction since the number of columns of
$T_{a}$ is finite. It means that the sequence $(T_{n})$ is finite.

\noindent Finally $T_{u}$ the last tableau of this sequence is necessarily a
column such that $\widehat{T}_{u}=T_{u}$ that is, $T_{u}$ is a column of
weight $0.$
\end{proof}

\begin{example}
The cocyclage operation is not authorized for $T=$%
\begin{tabular}
[c]{|l|l}\hline
$\mathtt{\bar{3}}$ & \multicolumn{1}{|l|}{$\mathtt{\bar{3}}$}\\\hline
$\mathtt{\bar{2}}$ & \multicolumn{1}{|l|}{$\mathtt{\bar{1}}$}\\\hline
$\mathtt{1}$ & \\\cline{1-1}%
\end{tabular}
. We have $\widehat{T}=$%
\begin{tabular}
[c]{|l|l}\hline
$\mathtt{\bar{3}}$ & \multicolumn{1}{|l|}{$\mathtt{\bar{2}}$}\\\hline
$\mathtt{2}$ & \\\cline{1-1}%
\end{tabular}
. Then $T_{1}=$%
\begin{tabular}
[c]{|l|l}\hline
$\mathtt{\bar{3}}$ & \multicolumn{1}{|l|}{$\mathtt{2}$}\\\hline
$\mathtt{\bar{2}}$ & \\\cline{1-1}%
\end{tabular}
and $T_{2}=$%
\begin{tabular}
[c]{|l|}\hline
$\mathtt{\bar{3}}$\\\hline
$\mathtt{3}$\\\hline
\end{tabular}
.
\end{example}

\subsection{Embedding of cyclage graphs}

Each connected component $\Gamma(T)$ contains a tableau $Y$ which admits no
cocyclage.\ This tableau is necessarily unique. Suppose that $Y\in
\mathbf{ST(}m\mathbf{)}$ .Then $\Gamma(T)\subset\mathbf{ST(}m\mathbf{)}$ thus
is finite. Moreover for any $Z\in\Gamma(T)$ there exits an integer $k$ such
that $U^{(k)}(Z)=T.$ This means that $\Gamma(T)$ has a tree structure.

\begin{proposition}
\label{prop_embedd}Let $\mu\in P_{n}$ and consider $T_{\mu}\in\mathbf{ST(}%
\mu\mathbf{)}$.\ Suppose that there exists $j\leq i\leq n$ such that
$\mu_{\overline{i}}>\mu_{\overline{j}}\geq0$.\ Set $\nu\in P_{n}$ defined by
$\nu_{\overline{k}}=\mu_{\overline{k}}$ for $k\neq i,j,$ $\nu_{\overline{i}%
}=\mu_{\overline{i}}-1$ and $\nu_{\overline{j}}=\mu_{\overline{j}}+1.$ Then
there exists a tableau $T_{\nu}\in\mathbf{ST(}\nu\mathbf{)}$ and a unique
embedding from $\Gamma(T_{\mu})$ to $\Gamma(T_{\nu})$ which commutes with $U$
and preserves the shape of the tableaux.
\end{proposition}

\begin{proof}
We have seen that $\Gamma(T_{\mu})$ has a finite number of vertices.\ Write
$m\geq n$ for the lowest integer such that $\Gamma(T_{\mu})$ is contained in
$\mathbf{ST(}m\mathbf{).\;}$By abuse of notation we also denote $\mu$ and
$\nu$ the weights of $P_{m}$ defined by $\mu=(0,...,0,\mu_{\overline{n}%
},...,\mu_{\overline{1}})$ and $\nu=(0,...,0,\nu_{\overline{n}},...,\nu
_{\overline{1}}).\;$Let $\sigma\in W_{m}$ such that $\sigma(\overline
{i})=\overline{m}$, $\sigma(\overline{j})=\overline{m-1}$ and $\sigma
(\overline{k})=\overline{k}$ for $k\neq i,j.$ Set $\sigma(\mu)=\mu^{\prime}$
and $\sigma(v)=v^{\prime}.$ Write $T_{\mu^{\prime}}=\sigma(T_{\mu}).\;$We have
$\Gamma(T_{\mu^{\prime}})\subset\mathbf{ST(}m\mathbf{).\;}$Then for any
$T\in\Gamma(T_{\mu^{\prime}}),$ $\widetilde{f}_{m-1}(\mathrm{w}(T))\neq0.$
Indeed the crystal graph $G_{m}^{(m-1)}$ obtained by erasing in $G_{m}$ all
the arrows of color $i\neq m-1$ and all the letters $x\notin\{\overline
{m},\overline{m-1},m-1,m\}$ is a $U_{q}(sl_{2})_{m-1}$-crystal where
$U_{q}(sl_{2})_{m-1}$ is the sub-algebra of $U_{q}(sp_{2m})$ isomorphic to
$U_{q}(sl_{2})$ generated by $e_{m-1},f_{m-1}$ and $t_{m-1}.$ The vertex
$\mathrm{w}(T)_{m-1}$ of $G_{m}^{(m-1)}$ obtained from $\mathrm{w}(T)$ is of
weight $(\mu_{\overline{m}}^{\prime},\mu_{\overline{m-1}}^{\prime})\neq
0$.\ Since $\mu_{\overline{i}}>\mu_{\overline{j}}\geq0$ we have $\mu
_{\overline{m}}^{\prime}>\mu_{\overline{m-1}}^{\prime}\geq0.$ Thus
$\mathrm{w}(T)_{m-1}$ is a highest weight vertex and there is an arrow of
color $m-1$ and length $\mu_{\overline{m}}^{\prime}-\mu_{\overline{m-1}%
}^{\prime}$ which starts from $\mathrm{w}(T).$

\noindent Now consider $T\in\Gamma(T_{\mu^{\prime}})$ such that
$U(T)=T^{\prime}$ is defined$.$ Write $\mathrm{w}(T)=x_{\ast}\mathrm{w}%
(T_{\ast}).\;$We must have $x_{\ast}\neq\overline{m}$ since the cocyclage
operation is authorized for $T$.\ Moreover $x_{\ast}\neq m.\;$Otherwise the
first column of $T^{\prime}$ would contain the letters $m$ and $\overline{m}$
(because $\mu_{\overline{m}}^{\prime}>0)$ and $T^{\prime}\notin\mathbf{ST(}%
m\mathbf{).}$ We are going to prove that
\begin{equation}
\widetilde{f}_{m-1}(x_{\ast}\mathrm{w}(T_{\ast}))=x_{\ast}\widetilde{f}%
_{m-1}(\mathrm{w}(T_{\ast})). \label{fi_on_w(tstar)}%
\end{equation}
It suffices to establish (\ref{fi_on_w(tstar)}) for $x_{\ast}=m-1.\;$When
$x_{\ast}=m-1$ there is no letter $\overline{m-1}$ in the second row of $T$
since $T^{\prime}\in\mathbf{ST(}m\mathbf{).}$ Thus all the letters
$\overline{m-1}$ of $T$ belong to its first row$.$ Denote by $T_{1}$ the sub
tableau of $T$ containing all the columns whose the lowest letter is
$\overline{m}.\;$The tableau $T$ can be regarded as the juxtaposition
$T_{1}T_{2}$ of the tableaux $T_{1}$ and $T_{2}$ where $T_{2}$ is the
sub-tableau obtained by considering the columns of $T$ which do not occur in
$T_{1}.$ Then $T_{1}$ do not contain any letter $m$ or $\overline{m-1}$ and
$T_{2}$ do not contain any letter $\overline{m}$. Suppose that $\widetilde
{f}_{m-1}(x_{\ast}\mathrm{w}(T_{\ast}))=\widetilde{f}_{m-1}(x_{\ast
})\mathrm{w}(T_{\ast}).$ Then with the notation of Note \ref{+-} we can write
$\rho(\mathrm{w}(T_{2}))=(+)^{s}$ since $x_{\ast}=m-1=+$ is not ignored during
the encoding procedure. The pairs $(+-)$ ignored are pairs $(m-1,\overline
{m-1})$ or $(m-1,m)$ for $\overline{m}$ do not belong to $\mathrm{w}%
(T_{2}).\;$Since $\rho(\mathrm{w}(T_{2}))$ contains only symbols $+,$ all the
letters $\overline{m-1}$ can be paired with letters $m-1.$ Thus the number of
letters $m-1$ in $\mathrm{w}(T_{2})$ is strictly greater than that of letters
$\overline{m-1}.$ It is also true for $\mathrm{w}(T)$ because $\mathrm{w}%
(T_{1})$ does not contain any letter $\overline{m-1}$.\ This contradicts the
inequality $\mu_{m-1}^{\prime}\geq0.$ Thus (\ref{fi_on_w(tstar)}) is true.

\noindent Denote by $V$ the symplectic tableau of reading $\widetilde{f}%
_{m-1}(\mathrm{w}(T_{\mu^{\prime}})).$ We are going to prove that $\Psi
:\Gamma(T_{\mu^{\prime}})\rightarrow\Gamma(V)$ defined by $\Psi(T)=S$ if and
only if $\mathrm{w}(S)=\widetilde{f}_{m-1}(\mathrm{w}(T))$ is an embedding
which commutes with $U$ and preserves the shape of the tableaux. We have
$\Psi(U(T))=P\left(  \widetilde{f}_{m-1}(\mathrm{w}(T_{\ast})x_{\ast})\right)
$. Suppose that there exits $p\in\{1,...,m\}$ such that $\nu_{\overline{p}%
}^{\prime}$ is equal to $r$ the number of columns of $S.$ Then since $U(T)$ is
defined and $\nu_{\overline{k}}^{\prime}=\mu_{\overline{k}}^{\prime}$ for any
$k\neq m,m-1$ we must have $p\in\{m,m-1\}.$ If $\nu_{\overline{m}}^{\prime}=r$
then we have $\mu_{\overline{m}}^{\prime}=r+1$ which is impossible for $T$
contains only $r$ columns.\ If $\nu_{\overline{m-1}}^{\prime}=r$ then we
obtain $\mu_{\overline{m}}^{\prime}\geq r$ since $\mu_{\overline{m-1}}%
^{\prime}=r-1$and $\mu_{\overline{m}}^{\prime}>\mu_{\overline{m-1}}^{\prime}%
$.\ This contradicts the fact that the cocyclage operation is authorized for
$T.$ Hence the cocyclage operation is authorized for $S$ and we can write%
\[
U(S)=U(\Psi(T))=P\left(  \xi(\widetilde{f}_{m-1}(x_{\ast}\mathrm{w}(T_{\ast
}))\right)  =P\left(  \xi(x_{\ast}(\widetilde{f}_{m-1}(\mathrm{w}(T_{\ast
})))\right)  =P\left(  \widetilde{f}_{m-1}(\mathrm{w}(T_{\ast}))x_{\ast
}\right)  .
\]
Thus we have to show that
\[
\widetilde{f}_{m-1}(\mathrm{w}(T_{\ast})x_{\ast})=\widetilde{f}_{m-1}%
(\mathrm{w}(T_{\ast}))x_{\ast}.
\]
By (\ref{TENS1}) it is equivalent to%
\begin{equation}
\varphi_{m-1}(\mathrm{w}(T_{\ast}))>\varepsilon_{m-1}(x_{\ast}). \label{ineg}%
\end{equation}
We have seen that $x_{\ast}\neq m$ and for $x_{\ast}\neq\overline{m-1}$
(\ref{ineg}) is true since $\varepsilon_{m-1}(x_{\ast})=0$ and $\varphi
_{m-1}(\mathrm{w}(T_{\ast}))\geq1.$ Suppose that $x_{\ast}=\overline{m-1}.$
Then the vertex of $G_{m}^{(m-1)}$ obtained from $\mathrm{w}(T_{\ast})$ as
above is of weight $(\mu_{\overline{m}},\mu_{\overline{m-1}}-1).$ Thus
$\varphi_{m-1}(\mathrm{w}(T_{\ast}))\geq\mu_{\overline{m}}-\mu_{\overline
{m-1}}+1\geq2$ for $\mu_{\overline{m}}>\mu_{\overline{m-1}}.$ So (\ref{ineg})
is satisfied. It is clear that $T$ and $\Psi(T)$ have the same
shape.\ Moreover by $\mathrm{(iii)}$ of Lemma \ref{lem_shape} $\Psi$ is the
unique map from $\Psi$ $\Gamma(T_{\mu^{\prime}})$ to $\Gamma(V)$ which
commutes with $U$ and preserves the shape of the tableaux. Finally, using
$\sigma^{-1}$ we obtain from $\Psi$ a unique embedding $\Psi_{\sigma}$
satisfying $\Psi_{\sigma}(T)=\sigma^{-1}\Psi\sigma(T)$ from $\Gamma(T_{\mu
})\ $to $\Gamma(T_{\nu})$ with $T_{\nu}=\sigma^{-1}\Psi\sigma(T_{\mu}).$
\end{proof}

\noindent Note that $T_{\nu}$ is not unique in general since $\Gamma(\nu)$ may
contain fewer connected components isomorphic to $\Gamma(T_{\nu})$.

\begin{corollary}
Let $\mu\in P_{n}^{+}$ and $T_{\mu}\in\mathbf{ST(}\mu\mathbf{).\;}$Write
respectively $m$ and $p$ for the sum of the non zero parts and the number of
zero parts in $\mu$.\ Define $\kappa\in P_{m+p}^{+}$ by $\kappa_{\overline{i}%
}=1$ for $p+1\leq i\leq m+p$ and $\kappa_{\overline{i}}=0$ otherwise. Then
there exists a tableau $T_{\kappa}\in\mathbf{ST(}\kappa\mathbf{)}$ and a
unique embedding of $\Gamma(T_{\mu})$ into $\Gamma(T_{\kappa})$ which commutes
with $U$ and preserves the shape of the tableaux.
\end{corollary}

\begin{proof}
The corollary directly follows by composing embeddings obtained in the
previous Proposition.
\end{proof}

\begin{example}
\label{big_graph}If $n=3$ and $\mu=(2,1,0)$ then $\kappa=(1,1,1,0).\;$The
cyclage graph $\Gamma\left(
\begin{tabular}
[c]{|l|l|l|l|l|}\hline
$\mathtt{\bar{3}}$ & $\mathtt{\bar{3}}$ & $\mathtt{\bar{2}}$ & $\mathtt{\bar
{1}}$ & $\mathtt{1}$\\\hline
\end{tabular}
\right)  $ of (\ref{graph(2,1,0)}) may be uniquely embedded in $\Gamma\left(
\begin{tabular}
[c]{|l|l|l|l|l|}\hline
$\mathtt{\bar{4}}$ & $\mathtt{\bar{3}}$ & $\mathtt{\bar{2}}$ & $\mathtt{\bar
{1}}$ & $\mathtt{1}$\\\hline
\end{tabular}
\right)  .$%

\[%
\begin{tabular}
[c]{ccccccccc}%
&  & $\vspace{0.1cm}%
\begin{tabular}
[c]{|l|}\hline
$\mathtt{\bar{4}}$\\\hline
$\mathtt{\bar{3}}$\\\hline
$\mathtt{\bar{2}}$\\\hline
$\mathtt{\bar{1}}$\\\hline
$\mathtt{1}$\\\hline
\end{tabular}
$ &  &  &  &  &  & \\
&  & $\vspace{0.1cm}\uparrow$ &  &  &  &  &  & \\
&  & $\vspace{0.1cm}%
\begin{tabular}
[c]{|l|l}\hline
$\mathtt{\bar{4}}$ & \multicolumn{1}{|l|}{$\mathtt{1}$}\\\hline
$\mathtt{\bar{3}}$ & \\\cline{1-1}%
$\mathtt{\bar{2}}$ & \\\cline{1-1}%
$\mathtt{\bar{1}}$ & \\\cline{1-1}%
\end{tabular}
$ &  &  &  &  &  & \\
& $\nearrow$ & $\vspace{0.1cm}\uparrow$ &  &  &  &  &  & \\%
\begin{tabular}
[c]{|l|l}\hline
$\mathtt{\bar{4}}$ & \multicolumn{1}{|l|}{$\mathtt{\bar{1}}$}\\\hline
$\mathtt{\bar{3}}$ & \multicolumn{1}{|l|}{$\mathtt{1}$}\\\hline
$\mathtt{\bar{2}}$ & \\\cline{1-1}%
\end{tabular}
&  & $\vspace{0.1cm}%
\begin{tabular}
[c]{|l|l}\hline
$\mathtt{\bar{4}}$ & \multicolumn{1}{|l|}{$\mathtt{\bar{1}}$}\\\hline
$\mathtt{\bar{3}}$ & \\\cline{1-1}%
$\mathtt{\bar{2}}$ & \\\cline{1-1}%
$\mathtt{1}$ & \\\cline{1-1}%
\end{tabular}
$ &  &  &  &  &  & \\
$\uparrow$ &  & $\vspace{0.1cm}\uparrow$ & $\nwarrow$ &  &  &  &  & \\%
\begin{tabular}
[c]{|l|l}\hline
$\mathtt{\bar{4}}$ & \multicolumn{1}{|l|}{$\mathtt{\bar{2}}$}\\\hline
$\mathtt{\bar{3}}$ & \multicolumn{1}{|l|}{$\mathtt{\bar{1}}$}\\\hline
$\mathtt{1}$ & \\\cline{1-1}%
\end{tabular}
&  & $\vspace{0.1cm}%
\begin{tabular}
[c]{|l|ll}\hline
$\mathtt{\bar{4}}$ & $\mathtt{\bar{1}}$ & \multicolumn{1}{|l|}{$\mathtt{1}$%
}\\\hline
$\mathtt{\bar{3}}$ &  & \\\cline{1-1}%
$\mathtt{\bar{2}}$ &  & \\\cline{1-1}%
\end{tabular}
$ &  &
\begin{tabular}
[c]{|l|l}\hline
$\mathtt{\bar{4}}$ & \multicolumn{1}{|l|}{$\mathtt{\bar{2}}$}\\\hline
$\mathtt{\bar{3}}$ & \\\cline{1-1}%
$\mathtt{\bar{1}}$ & \\\cline{1-1}%
$\mathtt{1}$ & \\\cline{1-1}%
\end{tabular}
&  &  &  & \\
$\uparrow$ &  & $\vspace{0.1cm}\uparrow$ &  & $\uparrow$ & $\nwarrow$ &  &  &
\\%
\begin{tabular}
[c]{|l|l|l}\hline
$\mathtt{\bar{4}}$ & $\mathtt{\bar{2}}$ & \multicolumn{1}{|l|}{$\mathtt{1}$%
}\\\hline
$\mathtt{\bar{3}}$ & $\mathtt{\bar{1}}$ & \\\cline{1-2}%
\end{tabular}
&  & $\vspace{0.1cm}$%
\begin{tabular}
[c]{|l|l}\hline
$\mathtt{\bar{4}}$ & \multicolumn{1}{|l|}{$\mathtt{\bar{2}}$}\\\hline
$\mathtt{\bar{3}}$ & \multicolumn{1}{|l|}{$\mathtt{1}$}\\\hline
$\mathtt{\bar{1}}$ & \\\cline{1-1}%
\end{tabular}
&  &
\begin{tabular}
[c]{|l|l}\hline
$\mathtt{\bar{4}}$ & \multicolumn{1}{|l|}{$\mathtt{\bar{3}}$}\\\hline
$\mathtt{\bar{2}}$ & \\\cline{1-1}%
$\mathtt{\bar{1}}$ & \\\cline{1-1}%
$\mathtt{1}$ & \\\cline{1-1}%
\end{tabular}
&  &
\begin{tabular}
[c]{|l|ll}\hline
$\mathtt{\bar{4}}$ & $\mathtt{\bar{2}}$ & \multicolumn{1}{|l|}{$\mathtt{1}$%
}\\\hline
$\mathtt{\bar{3}}$ &  & \\\cline{1-1}%
$\mathtt{\bar{1}}$ &  & \\\cline{1-1}%
\end{tabular}
&  & \\
& $\nearrow$ & $\vspace{0.1cm}\uparrow$ &  & $\uparrow$ &  & $\uparrow$ &  &
\\%
\begin{tabular}
[c]{|l|l|l}\hline
$\mathtt{\bar{4}}$ & $\mathtt{\bar{2}}$ & \multicolumn{1}{|l|}{$\mathtt{\bar
{1}}$}\\\hline
$\mathtt{\bar{3}}$ & $\mathtt{1}$ & \\\cline{1-2}%
\end{tabular}
&  & $\vspace{0.1cm}%
\begin{tabular}
[c]{|l|ll}\hline
$\mathtt{\bar{4}}$ & $\mathtt{\bar{2}}$ & \multicolumn{1}{|l|}{$\mathtt{\bar
{1}}$}\\\hline
$\mathtt{\bar{3}}$ &  & \\\cline{1-1}%
$\mathtt{1}$ &  & \\\cline{1-1}%
\end{tabular}
$ &  & $%
\begin{tabular}
[c]{|l|ll}\hline
$\mathtt{\bar{4}}$ & $\mathtt{\bar{3}}$ & \multicolumn{1}{|l|}{$\mathtt{1}$%
}\\\hline
$\mathtt{\bar{2}}$ &  & \\\cline{1-1}%
$\mathtt{\bar{1}}$ &  & \\\cline{1-1}%
\end{tabular}
$ &  &
\begin{tabular}
[c]{|l|l}\hline
$\mathtt{\bar{4}}$ & \multicolumn{1}{|l|}{$\mathtt{\bar{3}}$}\\\hline
$\mathtt{\bar{2}}$ & \multicolumn{1}{|l|}{$\mathtt{1}$}\\\hline
$\mathtt{\bar{1}}$ & \\\cline{1-1}%
\end{tabular}
&  & \\
& $\nearrow$ & $\vspace{0.1cm}\uparrow$ &  &  &  & $\uparrow$ & $\nwarrow$ &
\\%
\begin{tabular}
[c]{|l|lll}\hline
$\mathtt{\bar{4}}$ & $\mathtt{\bar{2}}$ & \multicolumn{1}{|l}{$\mathtt{\bar
{1}}$} & \multicolumn{1}{|l|}{$\mathtt{1}$}\\\hline
$\mathtt{\bar{3}}$ &  &  & \\\cline{1-1}%
\end{tabular}
&  & $\vspace{0.1cm}$%
\begin{tabular}
[c]{|l|l}\hline
$\mathtt{\bar{4}}$ & \multicolumn{1}{|l|}{$\mathtt{\bar{3}}$}\\\hline
$\mathtt{\bar{2}}$ & \multicolumn{1}{|l|}{$\mathtt{\bar{1}}$}\\\hline
$\mathtt{1}$ & \\\cline{1-1}%
\end{tabular}
&  &  &  &
\begin{tabular}
[c]{|l|ll}\hline
$\mathtt{\bar{4}}$ & $\mathtt{\bar{3}}$ & \multicolumn{1}{|l|}{$\mathtt{\bar
{1}}$}\\\hline
$\mathtt{\bar{2}}$ &  & \\\cline{1-1}%
$\mathtt{1}$ &  & \\\cline{1-1}%
\end{tabular}
&  &
\begin{tabular}
[c]{|l|l|l}\hline
$\mathtt{\bar{4}}$ & $\mathtt{\bar{3}}$ & \multicolumn{1}{|l|}{$\mathtt{\bar
{1}}$}\\\hline
$\mathtt{\bar{2}}$ & $\mathtt{1}$ & \\\cline{1-2}%
\end{tabular}
\\
& $\nearrow$ & $\vspace{0.1cm}\uparrow$ &  &  &  & $\uparrow$ &  & \\%
\begin{tabular}
[c]{|l|ll}\hline
$\mathtt{\bar{4}}$ & $\mathtt{\bar{3}}$ & \multicolumn{1}{|l|}{$\mathtt{\bar
{2}}$}\\\hline
$\mathtt{\bar{1}}$ &  & \\\cline{1-1}%
$\mathtt{1}$ &  & \\\cline{1-1}%
\end{tabular}
&  & $\vspace{0.1cm}%
\begin{tabular}
[c]{|l|l|l}\hline
$\mathtt{\bar{4}}$ & $\mathtt{\bar{3}}$ & \multicolumn{1}{|l|}{$\mathtt{1}$%
}\\\hline
$\mathtt{\bar{2}}$ & $\mathtt{\bar{1}}$ & \\\cline{1-2}%
\end{tabular}
$ &  &  &  &
\begin{tabular}
[c]{|l|lll}\hline
$\mathtt{\bar{4}}$ & $\mathtt{\bar{3}}$ & \multicolumn{1}{|l}{$\mathtt{\bar
{1}}$} & \multicolumn{1}{|l|}{$\mathtt{1}$}\\\hline
$\mathtt{\bar{2}}$ &  &  & \\\cline{1-1}%
\end{tabular}
&  & \\
$\uparrow$ &  & $\vspace{0.1cm}\uparrow$ &  &  &  &  &  & \\
$%
\begin{tabular}
[c]{|l|lll}\hline
$\mathtt{\bar{4}}$ & $\mathtt{\bar{3}}$ & \multicolumn{1}{|l}{$\mathtt{\bar
{2}}$} & \multicolumn{1}{|l|}{$\mathtt{1}$}\\\hline
$\mathtt{\bar{1}}$ &  &  & \\\cline{1-1}%
\end{tabular}
$ &  & $\vspace{0.1cm}%
\begin{tabular}
[c]{|l|l|l}\hline
$\mathtt{\bar{4}}$ & $\mathtt{\bar{3}}$ & \multicolumn{1}{|l|}{$\mathtt{\bar
{2}}$}\\\hline
$\mathtt{\bar{1}}$ & $\mathtt{1}$ & \\\cline{1-2}%
\end{tabular}
$ &  &  &  &  &  & \\
&  & $\vspace{0.1cm}\uparrow$ &  &  &  &  &  & \\
&  & $\vspace{0.1cm}%
\begin{tabular}
[c]{|l|lll}\hline
$\mathtt{\bar{4}}$ & $\mathtt{\bar{3}}$ & \multicolumn{1}{|l}{$\mathtt{\bar
{2}}$} & \multicolumn{1}{|l|}{$\mathtt{\bar{1}}$}\\\hline
$\mathtt{1}$ &  &  & \\\cline{1-1}%
\end{tabular}
$ &  &  &  &  &  & \\
&  & $\vspace{0.1cm}\uparrow$ &  &  &  &  &  & \\
&  &
\begin{tabular}
[c]{|l|l|l|l|l|}\hline
$\mathtt{\bar{4}}$ & $\mathtt{\bar{3}}$ & $\mathtt{\bar{2}}$ & $\mathtt{\bar
{1}}$ & $\mathtt{1}$\\\hline
\end{tabular}
&  &  &  &  &  &
\end{tabular}
\]
\end{example}

\section{A charge for symplectic tableaux $\label{sec_cha}$}

\subsection{Definition of $ch_{n}$}

\begin{definition}
\label{def_ch_col}Let $C$ be a column of weight $0.\;$Write $E_{C}%
=\{i\geq1,i\in C,i+1\notin C\}.\;$The charge $\mathrm{ch}_{n}(C)$ of the
column $C$ is%
\[
\mathrm{ch}_{n}(C)=2\sum_{i\in E_{C}}(n-i).
\]
\end{definition}

\noindent Note that for any column $T$ of weight $0$,
\[
\mathrm{ch}_{n+1}(t(C))=\mathrm{ch}_{n}(C)+2\mathrm{card}(E_{C})\text{ and
}\mathrm{ch}_{n+1}(t(C))=\mathrm{ch}_{n}(C).
\]
Moreover for any $i\geq1,$ we have $\varepsilon_{i}(\mathrm{w}((C)))=\left\{
\begin{tabular}
[c]{c}%
$1$ if $i\in E_{C}$\\
$0$ otherwise
\end{tabular}
\right.  $. Thus for any $n$-admissible column $C$ of weight $0$%
\begin{equation}
\mathrm{ch}_{n}(C)=2\sum_{i=1}^{n-1}(n-i)\varepsilon_{i}(\mathrm{w}(C)).
\label{chn_epsi}%
\end{equation}
Now consider $T\in\mathbf{ST}(\mu)$ with $\mu\in P_{n}^{+}$. Let
$\{T_{0},...,T_{p}\}$ with $T_{p}=C_{T}$ a column of weight $0$ be the above
sequence defined from $T.$

\begin{definition}
The charge $\mathrm{ch}_{n}(T)$ is%
\[
\mathrm{ch}_{n}(T)=\mathrm{ch}_{n}(C_{T})+p.
\]
\end{definition}

\noindent For any tableau $T$ we have by Lemma \ref{lem_shape}%
\[
\mathrm{ch}_{n+1}(t(T))=\mathrm{ch}_{n}(T).
\]

\subsection{Conjectures}

\begin{conjecture}
\label{conj1}Let $\Lambda_{k}^{A}$ and $\Lambda_{n-k}^{A}$ be respectively the
$k$-th and $(n-k)$-th fundamentals weights of $U_{q}(sl_{n})$.\ Set
$\lambda_{k}=\Lambda_{k}^{A}+\Lambda_{n-k}^{A}.\;$Then we have the equality:%
\[
K_{\Lambda_{2k},0}(q)=K_{\lambda_{k},0}^{A}(q^{2})
\]
where $K_{\lambda_{k},0}^{A}(q^{2})$ is the Kostka-Foulkes polynomial for the
root system $A_{n-1}$ corresponding to $\mu=0$ evaluated in $q^{2}.$
\end{conjecture}

\noindent Write $B(2k)_{0}=\{b\in B(\Lambda_{2k}),$ $\mathrm{wt}(b)=0\}$.\ By
identifying $U_{q}(sl_{n})$ with the subalgebra of $U_{q}(sp_{2n})$ generated
by the Chevalley's generators $e_{i},f_{i}$ and $t_{i},$ $i=1,...,n-1,$
$B(\Lambda_{2k})$ has a structure of crystal graph for $U_{q}(sl_{n})$
obtained by erasing all the arrows of color $0$ that we denote $B^{A}%
(\Lambda_{2k}).$ This graph decomposes into non isomorphic connected
components and in this decomposition $B_{0}(2k)$ is exactly the set of
vertices of weight $0$ of the connected component isomorphic to $B^{A}%
(\lambda_{k})$

\noindent By Theorem 5.1 of \cite{LLT} we can write
\[
K_{\lambda_{k},0}^{A}(q^{2})=\sum_{\mathrm{w}(C)\in B(2k)_{0}}q^{2\mathrm{d}%
^{\prime}(\mathrm{w}(C))}%
\]
with $\mathrm{d}^{\prime}(\mathrm{w}(C))=\sum_{i=1}^{n-1}(n-i)\varepsilon
_{i}(\mathrm{w}(C)).$ Then it follows from (\ref{chn_epsi}) that Conjecture
\ref{conj1} is equivalent to the equality%
\[
K_{\Lambda_{2k},0}(q)=\sum_{\mathrm{w}(C)\in B(2k)_{0}}q^{\mathrm{ch}_{n}%
(C)}.
\]
So by Corollary \ref{cor_h=2}, this conjecture is true for $k=1.$

\noindent More generally many computations suggest that $\mathrm{ch}_{n}$ is
an analogue for the root system $C_{n}$ of Lascoux-Sch\"{u}tzenberger's charge
on semi-standard tableaux.

\begin{conjecture}
\label{conj_charge}Consider $\lambda,\mu\in P_{n}^{+}.\;$Then%
\[
K_{\lambda,\mu}(q)=\sum_{\mathrm{w}(T)\in B(\lambda)_{\mu}}q^{\mathrm{ch}%
_{n}(T)}%
\]
where $B(\lambda)_{\mu}=\{T\in B(\lambda),\mathrm{wt}(T)=\mu\}.$
\end{conjecture}

\begin{example}
\ \ \ \ 

\begin{enumerate}
\item  Suppose $n=4,$ $\lambda^{(1)}=(2,1,1,1)$ and $\mu^{(1)}=(1,1,1,0).\;$%
There are $4$ tableaux in $\mathbf{ST}(4)$ of shape $\lambda^{(1)}$ and weight
$\mu^{(1)}$.\ They appear in the cyclage graph of Example \ref{big_graph}. Set
$C=\left(
\begin{tabular}
[c]{|l|}\hline
$\mathtt{\bar{4}}$\\\hline
$\mathtt{\bar{3}}$\\\hline
$\mathtt{\bar{2}}$\\\hline
$\mathtt{\bar{1}}$\\\hline
$\mathtt{1}$\\\hline
\end{tabular}
\right)  .\;$We have $\mathrm{ch}_{4}\left(  C\right)  =0$ since $\widehat
{C}=$%
\begin{tabular}
[c]{|l|}\hline
$\mathtt{\bar{4}}$\\\hline
$\mathtt{4}$\\\hline
\end{tabular}
. Hence the charges of these $4$ tableaux are $1,2,3$ and $4.\;$This gives the
Kostka-Foulkes polynomial $K_{\lambda^{(1)},\mu^{(1)}}(q)=q+q^{2}+q^{3}+q^{4}.$

\item  Suppose $n=3,$ $\lambda^{(2)}=(2,2,0)$ and $\mu^{(2)}=(0,0,0).\;$There
are $6$ tableaux in $\mathbf{ST}(3)$ of shape $\lambda^{(2)}$ and weight
$\mu^{(2)}$:%
\[
T_{1}=%
\begin{tabular}
[c]{|l|l|}\hline
$\mathtt{\bar{2}}$ & $\mathtt{1}$\\\hline
$\mathtt{\bar{1}}$ & $\mathtt{2}$\\\hline
\end{tabular}
,\text{ }T_{2}=%
\begin{tabular}
[c]{|l|l|}\hline
$\mathtt{\bar{2}}$ & $\mathtt{\bar{1}}$\\\hline
$\mathtt{1}$ & $\mathtt{2}$\\\hline
\end{tabular}
\text{, }T_{3}=%
\begin{tabular}
[c]{|l|l|}\hline
$\mathtt{\bar{3}}$ & $\mathtt{1}$\\\hline
$\mathtt{\bar{1}}$ & $\mathtt{3}$\\\hline
\end{tabular}
,\text{ }T_{4}=%
\begin{tabular}
[c]{|l|l|}\hline
$\mathtt{\bar{3}}$ & $\mathtt{\bar{1}}$\\\hline
$\mathtt{1}$ & $\mathtt{3}$\\\hline
\end{tabular}
,\text{ }T_{5}=%
\begin{tabular}
[c]{|l|l|}\hline
$\mathtt{\bar{3}}$ & $\mathtt{2}$\\\hline
$\mathtt{\bar{2}}$ & $\mathtt{3}$\\\hline
\end{tabular}
\text{ and }T_{6}=%
\begin{tabular}
[c]{|l|l|}\hline
$\mathtt{\bar{3}}$ & $\mathtt{\bar{2}}$\\\hline
$\mathtt{2}$ & $\mathtt{3}$\\\hline
\end{tabular}
.
\]
Note that $T_{2},T_{3}\in\Gamma\left(
\begin{tabular}
[c]{|l|l|l|l|}\hline
$\mathtt{\bar{1}}$ & $\mathtt{\bar{1}}$ & $\mathtt{1}$ & $\mathtt{1}$\\\hline
\end{tabular}
\right)  $ and $\Gamma(T_{1})$ is given in (\ref{graph(2,1,0)}).\ We obtain
$\mathrm{ch}_{3}\left(  T_{1}\right)  =2+\mathrm{ch}_{3}\left(
\begin{tabular}
[c]{|l|}\hline
$\mathtt{\bar{2}}$\\\hline
$\mathtt{\bar{1}}$\\\hline
$\mathtt{1}$\\\hline
$\mathtt{2}$\\\hline
\end{tabular}
\right)  =4,$ $\mathrm{ch}_{3}\left(  T_{2}\right)  =4+\mathrm{ch}_{3}\left(
\begin{tabular}
[c]{|l|}\hline
$\mathtt{\bar{3}}$\\\hline
$\mathtt{\bar{1}}$\\\hline
$\mathtt{1}$\\\hline
$\mathtt{3}$\\\hline
\end{tabular}
\right)  =8$ and $\mathrm{ch}_{3}\left(  T_{3}\right)  =6.$ Moreover
$T_{5}=t(T_{1})$ and $T_{6}=t(T_{2}).\;$Thus $\mathrm{ch}_{3}\left(
T_{5}\right)  =4-2=2$ and $\mathrm{ch}_{3}\left(  T_{6}\right)  =8-2\times
2=4.\;$By an easy computation we obtain $U^{(4)}(T_{4})=\left(
\begin{tabular}
[c]{|l|}\hline
$\mathtt{\bar{4}}$\\\hline
$\mathtt{\bar{1}}$\\\hline
$\mathtt{1}$\\\hline
$\mathtt{4}$\\\hline
\end{tabular}
\right)  .\;$Hence $\mathrm{ch}_{3}\left(  T_{4}\right)  =4+2(2-1)=6.$ This
gives the Kostka-Foulkes polynomial $K_{\lambda^{(2)},\mu^{(2)}}%
(q)=q^{2}+2q^{4}+2q^{6}+q^{8}.$
\end{enumerate}
\end{example}

\noindent\textbf{Remark:}

\noindent$\mathrm{(i):}$ Once $\mathrm{ch}_{n}$ defined on $\mathbf{ST,}$ it
is possible to define $\mathrm{ch}_{n}$ for any words of $\mathcal{C}_{\infty
}^{\ast}$ by setting
\[
\mathrm{ch}_{n}(w)=\mathrm{ch}_{n}(P(w)).
\]
Then given $w_{1},w_{2}\in\mathcal{C}_{n}^{\ast},$ the congruence $w_{1}%
\equiv_{n}w_{2}$ implies that $\mathrm{ch}_{n}(w_{1})=\mathrm{ch}_{n}(w_{2}),$
that is $\mathrm{ch}_{n}$ is a plactic invariant. We recover a property of the
Lascoux-Sch\"{u}tzenberger's charge $\mathrm{ch}_{A}$ for type $A$ \cite{LSc1}
\cite{Sc2}.\ Nevertheless, it seems difficult to define $\mathrm{ch}_{n}$
directly on words as it possible for $\mathrm{ch}_{A}.$ In \cite{LLT}, the
statistic $\mathrm{ch}_{A}$ is characterized in terms of the combinatorics of
crystal graphs.\ We have not found such a characterization for the symplectic
charge $\mathrm{ch}_{n}$.

\noindent$\mathrm{(ii):}$ It seems to be impossible to define a simple charge
statistic on $\mathbf{ST}(n)$ by using a cocyclage operation taking into
account the contraction relation (\ref{relcontra}) and relevant for computing
Kostka-Foulkes polynomials. Consider for example $T=%
\begin{tabular}
[c]{|l|l|l|l|}\hline
$\mathtt{\bar{1}}$ & $\mathtt{\bar{1}}$ & $\mathtt{1}$ & $\mathtt{1}$\\\hline
\end{tabular}
$ for $n=3.\;$If we apply cocyclages operations based on the complete
insertion scheme (with the contraction relations) we obtain the symplectic
tableaux of $\mathbf{ST}(3),$
\begin{tabular}
[c]{|l|ll}\hline
$\mathtt{\bar{1}}$ & $\mathtt{\bar{1}}$ & \multicolumn{1}{|l|}{$\mathtt{1}$%
}\\\hline
$\mathtt{1}$ &  & \\\cline{1-1}%
\end{tabular}
,
\begin{tabular}
[c]{|l|l|}\hline
$\mathtt{\bar{2}}$ & $\mathtt{\bar{1}}$\\\hline
$\mathtt{1}$ & $\mathtt{2}$\\\hline
\end{tabular}
,
\begin{tabular}
[c]{|l|ll}\hline
$\mathtt{\bar{2}}$ & $\mathtt{1}$ & \multicolumn{1}{|l|}{$\mathtt{2}$}\\\hline
$\mathtt{\bar{1}}$ &  & \\\cline{1-1}%
\end{tabular}
,
\begin{tabular}
[c]{|l|l}\hline
$\mathtt{\bar{2}}$ & \multicolumn{1}{|l|}{$\mathtt{1}$}\\\hline
$\mathtt{\bar{1}}$ & \\\cline{1-1}%
$\mathtt{2}$ & \\\cline{1-1}%
\end{tabular}
,
\begin{tabular}
[c]{|l|l}\hline
$\mathtt{\bar{3}}$ & \multicolumn{1}{|l|}{$\mathtt{3}$}\\\hline
$\mathtt{\bar{1}}$ & \\\cline{1-1}%
$\mathtt{1}$ & \\\cline{1-1}%
\end{tabular}
and
\begin{tabular}
[c]{|l|}\hline
$\mathtt{\bar{1}}$\\\hline
$\mathtt{1}$\\\hline
\end{tabular}
(since
\begin{tabular}
[c]{|l|}\hline
$\mathtt{\bar{3}}$\\\hline
$\mathtt{\bar{1}}$\\\hline
$\mathtt{1}$\\\hline
$\mathtt{3}$\\\hline
\end{tabular}
is not a $3$-admissible column). We know by Proposition \ref{prop_K_row} that
a charge for $T$ must necessarily be odd and by Corollary \ref{cor_h=2} a
charge for
\begin{tabular}
[c]{|l|}\hline
$\mathtt{\bar{1}}$\\\hline
$\mathtt{1}$\\\hline
\end{tabular}
must be even.\ So we can not deduce the charge of $T$ from that of
\begin{tabular}
[c]{|l|}\hline
$\mathtt{\bar{1}}$\\\hline
$\mathtt{1}$\\\hline
\end{tabular}
by simply counting the number of cocyclage operations.

\noindent$\mathrm{(iii):}$ Lascoux-Sch\"{u}tzenberger's proof of the equality%
\[
K_{\lambda,\mu}(q)=\sum_{\mathrm{w}(T)\in B(\lambda)_{\mu}}q^{\mathrm{ch}%
_{A}(T)}%
\]
for type $A$ is based on the Morris recurrence formula.\ We have seen that
Theorem \ref{Th_Morris} can be regarded as an analogue of this formula for
type $C_{n}.\;$It permits to decompose a Kostka-Foulkes polynomial for type
$C_{n}$ in terms of Kostka-Foulkes polynomials for type $C_{n-1}%
$.\ Unfortunately a charge statistic must take into account the contraction
relations to be compatible with the decomposition obtained in this way since
the partitions $\lambda$ such that $B(\lambda)$ appears in a decomposition\ of
type $B(\gamma)\otimes B((r)_{n-1})$ may be such that $\left|  \lambda\right|
<\left|  \mu\right|  .$ This is a reason why we are not able to deduce
Conjecture \ref{conj_charge} from Theorem \ref{Th_Morris}.

\bigskip

\bigskip
\end{document}